\documentclass[english, 12pt,twoside, a4paper]{article} 
\pdfoutput=1 
\RequirePackage[T1]{fontenc}
\RequirePackage[utf8]{inputenc}

\usepackage[
top=3cm,
bottom=2.8cm,
left=3cm,
right=3cm,
]{geometry}
\usepackage{setspace} %Réglage de l’interligne
\usepackage[nottoc, notlof, notlot]{tocbibind}
\usepackage{amsmath,amsfonts}
\usepackage{amsthm,amssymb}
\usepackage{letltxmacro}
\LetLtxMacro\amsproof\proof
\LetLtxMacro\amsendproof\endproof
\usepackage{thmtools}
\usepackage{empheq}
\allowdisplaybreaks
\usepackage[
unicode=true,
pdfstartview=FitV,
colorlinks=true,
citecolor=Leturquoise,
linkcolor=Leturquoise,
urlcolor=MFCB,
linktoc=page,
hyperindex=true,
pdfcreator={},
]{hyperref}
\usepackage{cleveref}
\crefname{figure}{Figure}{Figures}
\Crefname{figure}{Figure}{Figures}
\crefname{section}{Section}{Sections}
\crefname{CEx}{Counter-example}{Counter-examples}
\crefname{Ex}{Example}{Examples}
\crefname{Claim}{Claim}{Claims}
\crefname{Q}{Question}{Questions}

\usepackage{bm}% Pour le gras avec pdflatex
\usepackage{eulervm}
\usepackage{soul}
\sodef\allcapsspacing{\upshape}{0.15em}{0.65em}{0.6em}%
\sodef\lowsmallcapsspacing{\scshape}{0.075em}{0.5em}{0.6em}%   

\newcommand{\spacedlowsmallcaps}[1]{\MakeLowercase{\textsc{\lowsmallcapsspacing{#1}}}}
\usepackage[english]{babel}
\usepackage{adforn}
\usepackage{fancyhdr}
\usepackage[titles]{tocloft}
\usepackage{titlesec}
\titleformat{\section}
{\relax}{%
  \makebox[0pt][r]{\textbf{\Large \MakeLowercase{\scshape \thesection}}\hspace*{10pt}}%
}{1em}{\Large\bfseries \spacedlowsmallcaps}
\titleformat{\subsection}
{\relax}{\textbf{\large \thesubsection}}{1em}{\large \bfseries}
\titleformat{\subsubsection}
{\relax}{{{\bfseries \thesubsubsection}}}{1em}{\bfseries \normalsize}        
\titleformat{\paragraph}[runin]
{\normalfont\normalsize}{\theparagraph}{0pt}{\spacedlowsmallcaps}
\titlespacing*{\chapter}{0pt}{1\baselineskip}{1.2\baselineskip}
\titlespacing*{\section}{0pt}{1.25\baselineskip}{1\baselineskip} 
\titlespacing*{\subsection}{0pt}{1.25\baselineskip}{1\baselineskip}
\titlespacing*{\paragraph}{0pt}{1\baselineskip}{1\baselineskip}
\usepackage{enumitem}
\setlist{nosep}
\usepackage[x11names]{xcolor} %Pour changer les couleurs (sections, théorèmes...)
\usepackage{graphicx, subcaption} %Pour les includegraphics et compagnie
\usepackage{tikz} % Pour les schémas tikz
\usetikzlibrary{calc,math,patterns,patterns.meta}
\usetikzlibrary{arrows.meta} % Fioritures sur les flèches
\usepackage{tcolorbox}
\newtcolorbox{HypBox}{colframe=black, colback=black!0.3,boxrule=0.2mm, halign=left}
\AtBeginDocument{%
  \LetLtxMacro\proof\amsproof
  \LetLtxMacro\endproof\amsendproof
}
\declaretheorem[
style=definition,
thmbox={style=S,bodystyle=\normalfont \noindent},
name=Definition,
numberwithin=section,
refname={Definition,Definitions},
Refname={Definition,Definitions}]{Def}

\declaretheorem[
style=plain,
thmbox={style=S,bodystyle=\normalfont \noindent},
name=Theorem,
sibling=Def,
refname={Theorem,Theorems},
Refname={Theorm,Theorems}]{Th}

\declaretheorem[
style=plain,
thmbox={style=S,bodystyle=\normalfont \noindent},
name=Corollary,
sibling=Def,
sibling=Th,
refname={Corollary,Corollaries},
Refname={Corollary,Corollaries}]{Cor}

\declaretheorem[
style=plain,
thmbox={style=S,bodystyle=\normalfont \noindent},
name=Proposition,
sibling=Def,
refname={Proposition,Propositions},
Refname={Propositions,Propositions}]{Prop}

\declaretheorem[
style=plain,
thmbox={style=M,bodystyle=\normalfont \noindent},
name=Lemma,
sibling=Def,
heading={Lemma},
refname={Lemma,Lemmas},
Refname={Lemma,Lemmas}]{Lmm}
\renewcommand{\phi}{\varphi}
\renewcommand{\epsilon}{\varepsilon}
\theoremstyle{definition}
\newtheorem{Ex}[Def]{Example}
\newtheorem{CEx}[Def]{Counter-example}

\newtheorem{Q}[Def]{Question}
\newtheorem{Claim}[Def]{Claim}
\newtheorem{Rq}[Def]{Remark}

\newcommand{\bN}{\mathbb{N}}
\newcommand{\bZ}{\mathbb{Z}}
\newcommand{\bQ}{\mathbb{Q}}
\newcommand{\bR}{\mathbb{R}}

\newcommand{\bK}{\mathbb{K}}
\newcommand{\bF}{\mathbb{F}}

\newcommand{\calA}{\mathcal{A}}

\newcommand{\calC}{\mathcal{C}}
\newcommand{\calG}{\mathcal{G}}
\newcommand{\calH}{\mathcal{H}}

\newcommand{\calY}{\mathcal{Y}}
\newcommand{\mbfe}{\boldsymbol{e}}
\newcommand{\Tau}{T}%No \Tau
\newcommand{\Atlas}{\mathfrak{A}} %L’empreinte d’un sommet
\newcommand{\bdg}{\mathcal{X}} % L’immeuble
\newcommand{\myg}{X}% graphe étudié
\newcommand{\mY}{Y}% le Y R loc \myg
\newcommand{\CalY}{\mathcal{Y}}% le Y R loc l’immeuble
\newcommand{\Print}{\mathcal{P}} %L’empreinte d’un sommet
\newcommand{\PY}{P} %L’empreinte d’un sommet
\newcommand{\ratlas}{r_A} % rayon sur lequel l’atlas est bien défini
\newcommand{\rprint}{r_{\Print}} % Rayon sur lequel les print sont bien définies
\newcommand{\rX}{r_{\bdg}} % rayon sur lequel iota_y coincide avec Phi_y.
\newcommand{\rY}{r_{\mY}} % rayon sur lequel q{-1}iota kappa est une isom.
\newcommand{\Rbdg}{R_{\bdg}} % Caly et Rbdg localement l’immeuble
\newcommand{\Is}{\mathrm{Isom}} % Groupe d’isométrie
\newcommand{\Appart}{\mathcal{A}} % Appartement dans l’immeuble
\newcommand{\RInt}{\mathcal{O}} % Anneau des entiers
\newcommand{\mm}{\mathfrak{m}} % Ideal maximal
\definecolor{MarCB}{cmyk}{0,0.03,0.08,0.26} % BDB8AD 
\definecolor{MFCB}{cmyk}{0,0.06,0.20,0.6} % BDB8AD Mais en plus foncé
\colorlet{Leturquoise}{DeepSkyBlue4}
\colorlet{TurquoiseClair}{DeepSkyBlue4!20}
\colorlet{MyOrange}{DarkOrange3!85}
\tikzstyle{isometrie}=[->, >=latex, very thick]
\tikzstyle{qi}=[->, >=latex]
\tikzstyle{bqi}=[<->, >=latex]
\tikzstyle{nodeCzero}=[draw,circle, fill=blue!10]
\tikzstyle{nodeCun}=[draw,circle, fill=black!5]
\tikzstyle{unzero}=[dotted, thin]

\tikzstyle{vertex}=[draw,circle, fill=TurquoiseClair]
\tikzstyle{cercle}=[draw,circle]
\tikzstyle{bob}=[draw,circle, fill=TurquoiseClair]
\tikzstyle{bob1}=[draw,circle, fill=Leturquoise]
\tikzstyle{bob2}=[draw,circle, fill=MarCB]

\pgfmathsetmacro{\r}{1}
\pgfmathsetmacro{\h}{sqrt(3)/2}
\tikzstyle{geodesique}=[draw,thick,magenta!70]

\usepackage{titling}
\title{\texorpdfstring{Local-to-Global-rigidity of lattices in
    $SL_n(\bK)$}{Local-to-Global-rigidity of lattices in SLn(K)}}
\author{Amandine Escalier}
\date{\today}
\hypersetup{ %Les propriétés du PDF
pdfauthor={Amandine Escalier},
pdftitle={Local-to-Global-rigidity of lattices in SLn(K)},
pdfsubject={Local to Global behaviour of lattices in p-adic Lie groups},
pdfkeywords={groups, buildings, Bruhat, Tits, Kleiner, Leeb, isometries,
  rigidity, lattices, Lie, Cayley, Abels, graphs}
}

\begin{document}
\maketitle
\begin{abstract}
  A vertex-transitive graph $\calG$ is called \textit{Local-to-Global rigid} if there
  exists $R>0$ such that every other graph whose balls of radius $R$ are isometric
  to the balls of radius $R$ in $\calG$ is covered by $\calG$. An example of
  such a graph is given by the Bruhat-Tits building of $PSL_n(\bK)$ with
  $n\geq 4$ and $\bK$ a non-Archimedean local field of characteristic zero. In
  this paper we extend this rigidity property to a class of graphs 
  quasi-isometric to the building including torsion-free lattices of
  $SL_n(\bK)$.\\
  The proof is the occasion to prove a result on the local structure of the
  building. We show that if we fix a $PSL_n(\bK)$-orbit in it, then a vertex
  is uniquely determined by the neighbouring vertices in this orbit.
\end{abstract}

% \vspace{0.25cm}
% {\small \noindent\textbf{Classification} 20F65\\
% \textbf{Keywords} lattices, buildings, rigidity, local fields}
% \vspace{0.5cm}
\tableofcontents

% Introduction
\section{Introduction} 
A recurring theme in geometric group theory is that \emph{local} properties of an
object can have \emph{global} implication for its geometry. A classical example is
given by Lie groups and their locally defined Lie algebras. Another striking
illustration is provided by the work of Tits \cite{Tits} who 
gave a local characterization of a particular family of graphs called “buildings
of type $\tilde{A}_{d-1}$” (see \cref{BTBuilding} for a definition). 
Precisely, graphs and their 
local-to-global properties are the objects we focus on in this article. All
graphs will be equipped with the usual metric, fixing the length of an edge to one.

A natural local condition to impose on a graph is to be $d$-regular for
some $d\in \bN$, which means that all the vertices must have degree $d$. A
well-known result about such a graph is that the $d$-regular tree is its
universal convering. This is a first example of a global information
deduced only by a local knowledge of the graph.

One can now ask what happens if we impose a local condition which is stronger
than $d$-regularity. We formalize this in the next definition.

\begin{Def} Let $R>0$ and let $X$ and $Y$ be two graphs.\\
  We say that $Y$ is \textbf{R-locally} $X$ if every ball of radius $R$ in $Y$
  is isometric to a ball of radius $R$ in $X$.\\
  If $Y$ is $R$-locally $X$ and $X$ is $R$-locally $Y$ then we say that they are
  \textbf{R-locally the same.} 
\end{Def}

\begin{Ex}
  In the following example, $B_X(x_0,2)$ is isometric to $B_Y(y_0,2)$.
  \begin{figure}[htbp]
    \centering
    \begin{tikzpicture}[scale=0.9]
      \foreach \x in {0,30,...,360}{
        \draw (\x:2) -- (\x+30:2);
        \node[vertex] at (\x:2) {};
      };
      \foreach \x in {30,60,...,150}{
        % \draw (\x:2) -- (\x+30:2);
        \node[bob1] at (\x:2) {};
      };
      \node (X) at (-1.5,-2) {$X$};
      \node (x0) at (90:2.5) {$x_0$};
      \begin{scope}[shift={(4,0)}]
        \draw[dashed] (8,0) -- (9,0);
        \draw[dashed] (0,0) -- (-1,0);
        \foreach \x in {0,1,...,7}{
          \draw (\x,0) -- (\x+1,0);
          \node[vertex] at (\x,0) {};
        };
        \node[vertex] at (8,0) {};
        \foreach \x in {2,3,...,6}{
          \node[bob1] at (\x,0) {};
        };
        \node (Y) at (5,-2) {$Y$};
      \node (y0) at (4,0.5) {$y_0$};
      \end{scope}
    \end{tikzpicture}
    \caption{Two graphs 2-locally the same.}
    \label{fig:Rlocmeme}
  \end{figure}
\end{Ex}

The previous covering result on the $d$-regular tree is a first example of a more
general notion called the \textit{Local-to-Gobal rigidity}, also named
\textit{LG-rigidity}.

\begin{Def} Let $R >0$.
  We say that $X$ is \textbf{Local-to-Global-rigid\index{LG-rigidity} at scale
  $R$} (or \textbf{$R$-LG-rigid} for short) if every graph $Y$ which is
  $R$-locally $X$ is covered by $X$. \\
  We say that a graph $X$ is \textbf{LG-rigid} if there exists $R>0$ such
  that $X$ is $R$-LG-rigid.
\end{Def}

\begin{Ex} Benjamini and Ellis \cite{BELocLat} showed that for any $d \geq 2$ the Cayley graph of
  ${\bZ}^d$ endowed with its usual generating set is $3$-LG-rigid. They also
  proved that $3$ was optimal showing that $\bZ^3$ is not LG-rigid at scale $2$.
\end{Ex}

\begin{Ex}\label{Ex:QuasiArbre}
  De la Salle et Tessera
  \cite[Theorem C]{DlST2015} proved that every cocompact graph
  quasi-isometric to a tree is LG-rigid.
\end{Ex}
Benjamini \cite{B} and Georgakopoulos \cite{G} conjectured that any Cayley
graph of a finitely presented group is LG-rigid at some scale $R>0$.
That conjecture was proven to be false in \cite[Theorem B]{DlST2015}, where the
authors built counter-examples using groups with \emph{torsion} elements.
\begin{CEx} \label{CEx:torsion}
  The groups $F_2 \times F_2 \times \bZ/2\bZ$ and $SL_{4}(\bZ)$ admit Cayley
  graphs that are not LG-rigid.
\end{CEx}

Remark here that we do not state that every Cayley graph of these groups is
non-LG-rigid, but that each group \textit{admits} a non-LG-rigid 
Cayley graph. Indeed, in \cite[Theorem
J]{DlST2015} the authors also showed that 
every finitely presented group with an element of infinite order has a Cayley
graph which is LG-rigid. Hence, LG-rigidity
for a Cayley graph depends on the generating set. In particular LG-rigidity is
not invariant under quasi-isometries. 

With a little bit more of material, we will be able to give a topological
interpretation of Local-to-Global rigidity (see \cpageref{LGRtopo}).
\smallskip

That rigidity notion can be refined in what is called the \textit{Strong
  Local-to-Global rigidity}, also named \textit{SLG-rigidity}. 

\begin{Def} Let $r,R>0$. We say that $X$ is
  \textbf{SLG-rigid\index{SLG-rigidity} at scale $(r,R)$} if 
  for all $Y$ which is $R$-locally $X$ and for all isometry $f$ from
  $B_X(x,R)$ to $B_Y(y,R)$, the restriction of $f$ to $B_X(x,r)$ extends to a
  covering of $Y$ by $X$.\\
  We say that $X$ is \textbf{SLG-rigid} if there exist two radii $r$ and $R$
  such that $X$ is SLG-rigid at scale $(r,R)$.
\end{Def}

Such a refinement is far more than just a subtlety: it actually proves necessary to obtain our
main result (see \cpageref{firstmentionofSLGX} for more details).

The following proposition gives us many examples of SLG-rigid graphs.

\begin{Prop}[de la Salle, Tessera {\cite[Proposition
    3.8]{DlST2015}}] \label{Prop:LGRSLGR} A graph with cocompact 
  isometry group is LG-rigid if and only if it is SLG-rigid.
\end{Prop}

For example, any LG-rigid Cayley graph is actually SLG-rigid. In the same
article, de la Salle and Tessera proved a powerful condition 
relating to the isometry group of a Cayley graph. We will refer to the isometry
group of a Cayley graph $(\Gamma,S)$ as $\Is(\Gamma,S)$. 
\begin{Th}[de la Salle, Tessera {\cite[Theorem E]{DlST2015}}]\label{Th:DlSTdiscrete}
  Let $\Gamma$ be a finitely presented group and $S$ be a symmetric generating
  set and denote by $(\Gamma,S)$ the corresponding Cayley graph. If $\Is(\Gamma,S)$ is
  discrete, then $(\Gamma,S)$ is SLG-rigid.
\end{Th}

As stated in \cite[Corollary F]{DlST2015}, we can deduce two new
classes of examples from the above theorem. But before, let us introduce what we
call \emph{LG-rigid groups}.

\begin{Def}
  We say that a finitely presented group is \textbf{LG-rigid} (resp.
  \textbf{SLG-rigid}) if all its Cayley graphs are LG-rigid (resp.
  {SLG-rigid}). 
\end{Def}

\begin{Ex}
 Torsion-free groups of polynomial growth are SLG-rigid. 
\end{Ex}
\begin{Ex}
 Torsion-free, non-virtually free lattices in connected simple real Lie groups are SLG-rigid. 
\end{Ex}

So far, the graphs chosen as examples are mostly Cayley graphs, but these are not
the only LG-rigid ones. Indeed, besides the case of quasi-trees seen above, another
interesting example is given by \emph{Bruhat-Tits buildings} (see
\cref{BTBuilding} for a definition).

\begin{Th}[de la Salle, Tessera, {\cite[Theorem
    0.1]{DlST2016}}] \label{Th:LGRImmeuble}
  Let $\bK$ be a non-Archimedean local skew field.\\
  If $\bK$ has positive characteristic and $n\geq 3$, then the Bruhat-Tits building of
  $PSL_n(\bK)$ is not LG-rigid.\\
  If $\bK$ has characteristic zero and $n\geq 4$, then the Bruhat-Tits
  building of $PSL_n(\bK)$ is SLG-rigid.
\end{Th}

Keeping in mind the above theorem, consider the following question asked in \cite{DlST2015}. 
\begin{Q}\label{Q:lattices}
  \textit{Among lattices in semi-simple Lie groups, which ones are LG-rigid?}
\end{Q}
This question concerns \textit{real} Lie groups but one can also wonder what 
happens for the \textit{$p$-adic} case. Indeed, by a well known result of Svarc
and Milnor, any lattice of $SL_n(\bK)$ is quasi-isometric to the associated building
(see \cref{Lmm:GammaQIX} for more details). The fact that such a lattice is
“almost” a building encouraged us to study the $p$-adic version of
\cref{Q:lattices}.

\begin{Q}\label{Q:plattices}
  \textit{Among lattices in $p$-adic Lie groups, which ones are LG-rigid?}
\end{Q}
% --
De la Salle and Tessera showed \cite{DlST2016} that if $\bK$ has positive
characteristic, then there exist $p$-adic lattices that are torsion-free,
cocompact but not LG-rigid.
\begin{Ex} Let $n\geq 3$. There exists a torsion-free cocompact lattice in
  $PGL_n(\bF_p)$ that is not LG-rigid.
\end{Ex}
% --
When $\bK$ is a non-Archimedean local skew field of characteristic zero,  
an element of response to \cref{Q:plattices} is provided by our first result hereunder.
\begin{Th} \label{Th:LGRReseaux} Let $n \neq 3$ and $\bK$ be a non-Archimedean local skew
  field of characteristic zero.\\
  The torsion-free lattices of $SL_n(\bK)$ are SLG-rigid. 
\end{Th}

This result is actually a corollary of our main theorem below which
goes beyond the lattices framework and gives a rigidity result in a more general
case.   
\begin{Th}
  \label{Th:BigRigidity} Let $n\neq 3$ and $\bK$ be a non-Archimedean local skew
  field of characteristic zero. Let $\bdg$ be the
  Bruhat-Tits building of $PSL_n(\bK)$ and $\myg$ be a transitive graph. If
  $\myg$ verifies that
  \begin{itemize}
  \item There is an injective homomorphism $\rho$ from $\Is(\myg)$ to $\Is(\bdg)$ such that
    $\rho(\Is(\myg))$ is of finite index in $\Is(\bdg)$;
  \item There is a $\Is(\myg)$-equivariant injective quasi-isometry $q$ from $\myg$ to $\bdg$;
  \end{itemize}
  then $\myg$ is SLG-rigid.
\end{Th}

Let us discuss the hypothesis, starting with the assumption made on $n$. If
$n=2$ then $\bdg$ is the $(p+1)$-regular tree, thus by \cref{Ex:QuasiArbre} any
graph quasi-isometric to $\bdg$ is LG-rigid which proves the theorem. Now, as we
will see in the sketch of the proof, the main tool of our demonstration is the
LG-rigidity of the building. But if $n=3$ the question of the rigidity of $\bdg$
is still open. Indeed in that case a lot of flexibility seems to be allowed (see
\cite{Barre}). Thus our demonstration deals mainly with the case where $n\geq
4$.

Secondly, let us look at the hypothesis made on the characteristic of $\bK$.
According to \cite[Theorem 0.4]{DlST2016} and more precisely according to its
\emph{proof}, we get \cref{CEx:FpT} below. It implies in particular that if we omit
the characteristic zero hypothesis, then \cref{Th:BigRigidity,Th:LGRReseaux} are not true.
\begin{CEx}\label{CEx:FpT}
  There exists a non-LG-rigid torsion-free cocompact lattice in
  $PGL_n\big(\bF_p((T))\big)$.
\end{CEx}

Finally, before moving to the sketch of the proof let us discuss the hypothesis made
on the torsion in \cref{Th:LGRReseaux}. First, introducing torsion in a group is
in some case a useful way to build 
non-LG-rigid graphs. Indeed the \cref{CEx:torsion} is built this way.
Second, in order to link $(\Gamma,S)$ to $\bdg$ we will need an injection of
$\Is(\Gamma,S)$ into $\Is(\bdg)$. Using a famous result of Kleiner and Leeb we will
show that $\Is(\Gamma,S)$ acts on the buildings by isometries. The injection into
$\Is(\bdg)$ will then be allowed by the following
proposition.

\begin{Prop}[de la Salle, Tessera {\cite[Proposition 6.2]{DlST2015}}]
  Let $\Gamma$ be an infinite, torsion-free, finitely generated group and let
  $S$ be a finite symmetric generating subset of $\Gamma$. Then the isometry
  group of $(\Gamma,S)$ has no non-trivial compact normal subgroup.
  \label{DlSTProp62}
\end{Prop}
For more details on how we use this proposition, see the proof of \cref{Lmm:KLaffine}.  

\paragraph{Sketch of the proof of the theorem} As stated in the
discussion below \cref{Th:BigRigidity}, the proof deals mainly with the case
where $n\geq4$. So, Let $n\geq 4$ and $\bK$ be
non-Archimedean local skew field of characteristic zero 
and denote by $\bdg$ the Bruhat-Tits building of $PSL_n(\bK)$. Let $\myg$ be
the studied graph and $\mY$ be a graph $R$-locally the same as $\myg$ and denote
by $q$ a quasi-isometry from $\myg$ to $\bdg$. The main
idea of the proof is to use the rigidity of $\bdg$ to 
build the wanted covering from $\myg$ to $\mY$ (see \cref{fig:sketchproof}),
thus we need to build a graph locally the same as $\bdg$. We will denote such a graph $\CalY$.

\begin{figure}[htbp]
  \centering
    \begin{tikzpicture}[scale=1.3]
      % \myg and Y rloc myg
      \node[Leturquoise] (myg) at (0,2) {$\myg$};
      \node[Leturquoise] (mY) at (0,0) {$\mY$};
      \draw[<->, >=latex, dashed] (myg.south) -- (mY.north) node[left, midway] {$R$-loc};
      %Immmeuble and QI 
      \begin{scope}[shift={(2,0)}]
        \node[MFCB] (CalY) at (2,0) {$\calY$};
        \node[MFCB] (bdg) at (2,2) {$\bdg$};
        \draw[->, >=latex, very thick, black!50] (bdg.south west)to[bend right=90] (CalY.north west) ;
        \node[anchor=east, black!50] (cov) at (1.3,1) {\textbf{covering}};
      \end{scope}
      \draw[{Hooks[right]}->, >=latex] (myg.east) -- (bdg.west)
      node[above, anchor=south, midway] {$q$}
      node[below, anchor=north, midway] {quasi-isom.};
      \draw[<->, >=latex, dashed] (bdg.south) -- (CalY.north) node[right, midway] {$\Rbdg$-loc};
      \draw[->, >=latex, thick, dotted, black!50] (4,-0.7) -- (0,-0.7)%
      node[above, midway] {Goal:}
      node[below, midway] {induce a convering from $\myg$ to $\mY$};
    \end{tikzpicture}  
  \caption{Sketch of the proof}
  \label{fig:sketchproof}
\end{figure}

Moreover, for the rigidity of the
building to induce a covering between $\myg$ and $\mY$, we want $\CalY$
to contain a copy of the vertices of $\mY$. Hence the goal is to define the
vertices of $\CalY$
to be composed of the vertices of $\mY$ and a copy of each vertex in
$\bdg \backslash q(\myg)$
and define the edges to correspond to edges in $X$.
With such a description $\CalY$ is a “hybrid” graph and to define its edges we might
need to understand how to link a vertex coming from $\mY$ to a vertex coming
from $\bdg$. Hence, to avoid such a hybridation we chose to define the vertices
only with informations encoded in $\mY$. That is why we introduce
the notion of \emph{print} in the building (see
\cref{sec:localcharacterizationbdg}). It allows 
us to characterize a vertex in $\bdg$ by a set of neighbouring vertices in $im(q)$
and, using a well chosen set of isometries from $\mY$ to $\myg$, to transfer this print notion to
$\mY$. Each print in $\mY$ corresponds to a vertex in $\bdg \backslash
q(\myg)$. The vertices of the wanted graph $\CalY$ will be composed of the vertices of
$\mY$ and of prints in $\mY$. It will now be easier to build edges between these
vertices; the key argument to construct such edges is presented in \cref{subsec:Extofisom}.

Using the rigidity of the building we will obtain an isometry between $\bdg$ and
$\CalY$. To conclude the proof we will show that this isometry induces the
wanted covering between $\mY$ and $\myg$.

\paragraph{Organization of the paper} The first section is devoted to the
definition of our framework. We recall some material about Bruhat-Tits
buildings and large scale simple connectedness and present a fundamental
result on isometries’ extension.
The second and third sections are devoted to the proof of \cref{Th:BigRigidity}.
In the second section we develop the necessary engineering to build a graph
locally the same as the building —this is where we define and study prints—
while in the third one we use the rigidity of the building to prove the rigidity
of the studied graph.
We prove \cref{Th:LGRReseaux} in the fourth section where we check that the
lattice verifies the hypothesis of our main theorem.

\begin{center}
 \textbf{Acknowledgments} 
\end{center}
I would like to thank my advisor, Romain Tessera, under whose
supervision the work presented in this article was carried out. I thank him for suggesting the
topic and sharing his precious expertise during numerous discussions. I would also like
to thank Georges Skandalis for his useful advice, comments and corrections and Jean
Lécureux, Nicolas Radu and Sylvain Barré for helpful 
discussions on buildings. Finally, I thank Claire Debord for her support and
constructive remarks and the anonymous referee for fruitful questions and comments. 

% Section 1 : Cadre
\section{Framework} 
Let us start by setting up the framework of the next sections. We first recall
some material about Bruhat-Tits buildings, and large scale simple
connectedness. Then we present a useful tool concerning the extension of
isometries. We conclude by a result one step further in to the proof of our main
theorem, linking the $PSL_n(\bK)$-orbits in the building and the image
$q(\myg)$ of the graph studied.

\subsection{Bruhat-Tits building}
\label{BTBuilding} Let $n \geq2$. Since it is the object at the
center of our proof, let us recall the description 
of the Bruhat-Tits building associated to $PSL_n(\bK)$ where $n\geq 2$, see \cite{AB} for more
details.

\paragraph{Non-Archimedean local skew fields}
Let $\bK$ be a field (not necessarily commutative). A \emph{discrete valuation}
on $\bK$ is a surjective homomorphism $v \ : \ \bK^{*} \rightarrow \bZ$
satisfying $v(x+y) \geq \min\{v(x),v(y)\}$ for all $x,y \in \bK^{*}$ such that $x+y\neq 0$.
If $\bK$ is endowed with such a valuation, we can extend
$v$ on all $\bK$ by setting $v(0)=+ \infty$. We say that $\bK$ is a
\emph{non-Archimedean local skew field} if it is locally compact for the
topology associated to a discrete valuation.

\begin{Ex}
  If $\bK=\bQ$ and $p$ is a prime, then every $x\in \bK$ can be written as $x=p^n
  a/b$ where $a$ and $b$ are integers non-divisible by $p$. The map defined by
  $v(p^na/b):=n$ is a discrete valuation over $\bK$. The field $\bQ_p$ is the
  completion of $\bQ$ with respect to the \emph{p-adic absolute value} defined
  by $|x|_p=p^{-v(x)}$.
\end{Ex}
\begin{Ex}
  Let $\bK=\bF_p((T))$, the field of \emph{formal Laurent series} over $\bF_p$.
  Denote by $f=\sum_{k\in \bZ}a_k T^k$ an element in $\bF_p((T))$ then the map
  defined by $v(f):=\min\left\{ k \, : \, a_k \neq 0 \right\}$ is a valuation
  over $\bK$.    
\end{Ex}

Let $\RInt$ denote the \emph{ring of integers} of $\bK$ with respect to $v$,
that is to say $\RInt:=\{x \in \bK \ : \ v(x) \geq
0\}$. This ring has a unique \emph{prime ideal} $\mm:=\{x \in K \ : \ v(x)>0\}$.
Finally, let $\pi$ be a generator of $\mm$ as an $\RInt$-module.

\begin{Ex}
  If $\bK=\bQ_p$ then its ring of integers is $\RInt= \bZ_p$. Moreover $\mm= p\bZ_p$
  and $\pi =p$.
\end{Ex}
\begin{Ex}
  If $K=\bF_p((T))$ then $\RInt=\bF_p[[T]]$. Moreover $\mm=X \bF_p[[T]]$ and
  $\pi= X$.
\end{Ex}

\paragraph{Buildings}
Let $\bK$ be a non-Archimedean local skew field endowed with a valuation
$v$. An $\RInt$-lattice of $\bK^n$ is an $\RInt$-submodule which generates $\bK^n$
as a $\bK$ vector space. Such a lattice can be written as $\RInt e_1 + \cdots
+ \RInt e_n$ for a basis $(e_1, \ldots, e_n)$ of $\bK^n$. Since for any $a \in \bK^*$ and
any lattice $L$, the module $aL$ is also a lattice, we can define the
equivalence relation of \textit{lattices modulo homothety.} We denote by $[L]$
the class of a lattice $L$.\\
The Bruhat-Tits building of $PSL_n(\bK)$ is a simplicial complex of
dimension $n-1$ denoted by $\hat{\bdg}$ whose $1$-skeleton (denoted by $\bdg$) is described as
follows. The vertices are the classes of $\RInt$-lattices modulo homothety. Two
vertices $x_1$ and $x_2$ are linked by an edge if there exists representatives
$L_1$ of $x_1$ and $L_2$ of $x_2$ such that:
\begin{equation*}
 pL_1 \subset L_2 \subset L_1. 
\end{equation*}

\begin{figure}[htbp]
  \centering
  \begin{subfigure}[b]{0.42\textwidth}
    \includegraphics[height=\textwidth]{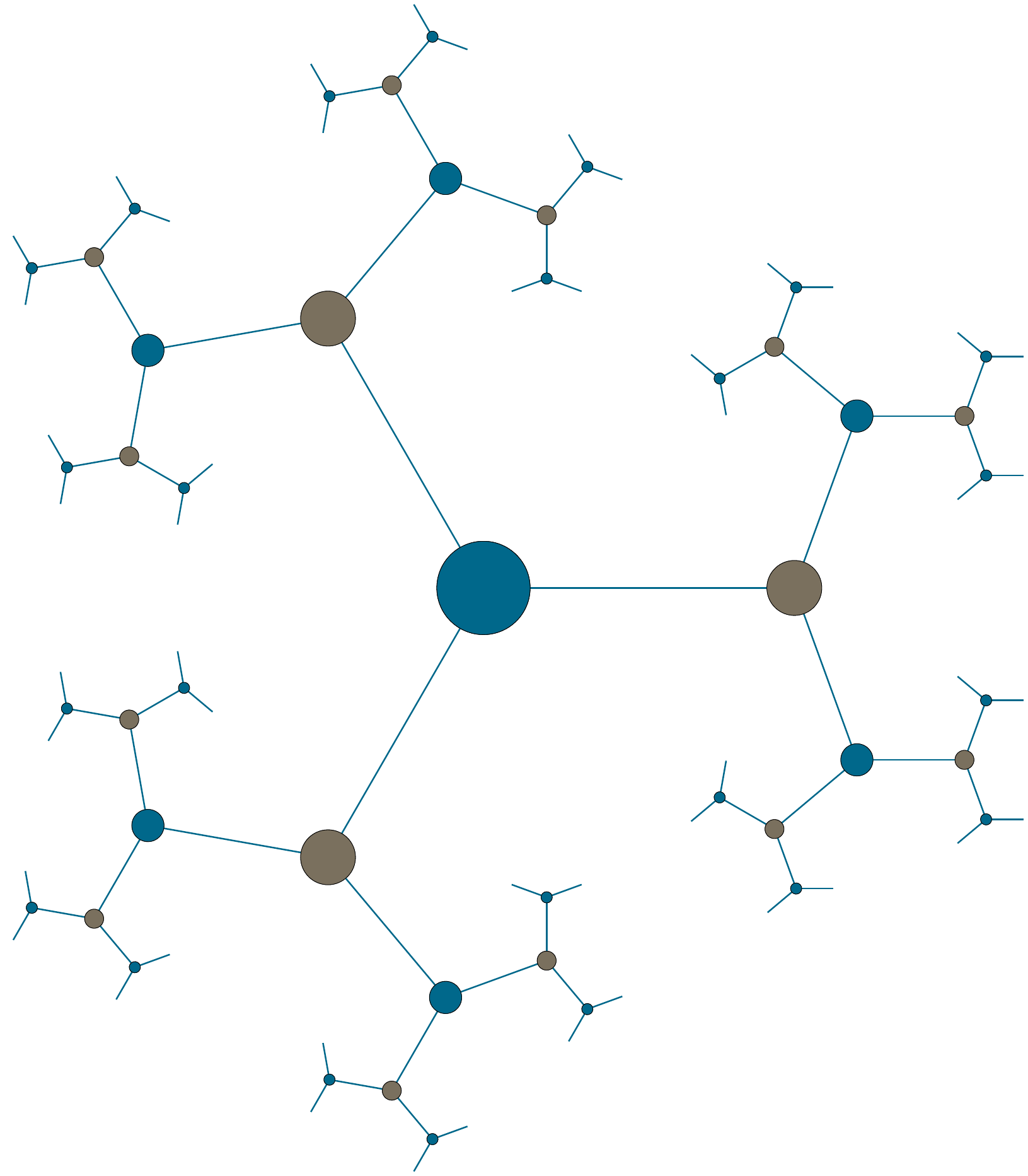}  
    \caption{The building has two $SL_2(\bQ_2)$-orbits}
    \label{fig:orbitsbuilding}
  \end{subfigure}
  \hfill
  \begin{subfigure}[b]{0.42\textwidth}
    \includegraphics[height=\textwidth]{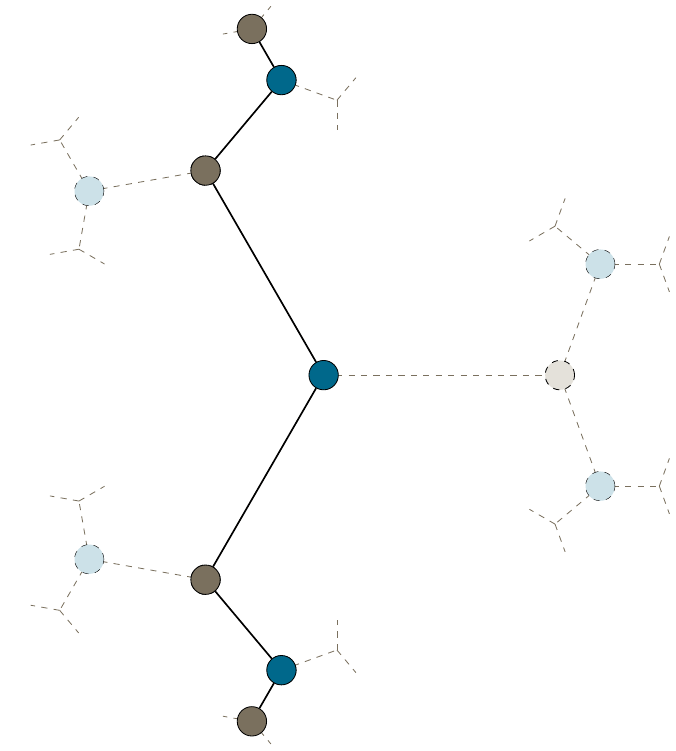}  
    \caption{Representation of one apartment}
    \label{fig:ApptArbre}
  \end{subfigure}
  \caption{The building of $PSL_2(\bQ_2)$}
\end{figure}

\begin{Ex}
 One can show that the building of $PSL_2(\bQ_p)$ is a $(p+1)$-regular tree.
 \Cref{fig:orbitsbuilding} gives a representation of the building when $p=2$. 
\end{Ex}

\paragraph{Orbits and types}
The usual action of $GL_n(\bK)$ on $\bK^n$ induces an action of
$PGL_n(\bK)$ on $\bdg$ by isometry. Since $GL_n(\bK)$ acts transitively on the
bases, the action of $PGL_n(\bK)$ on the vertices of $\bdg$ is also transitive.
\smallskip

\begin{figure}[htbp]
  \centering
  \includegraphics[width=0.6\textwidth]{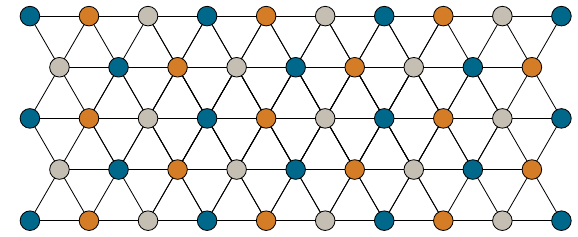}
  \caption{Apartment in the building of $PSL_3(\bQ_2)$. The colors correspond to
  $SL_3(\bQ_2)$-orbits.}
  \label{fig:ApptPSL3}
\end{figure}

If $L=\oplus_{i} \RInt e_i$ is a
lattice we define its \textit{type} to be $v\left(\det(e_1,\ldots, e_n)
\right)$. Since:
\begin{equation*}
  \forall a \in \bK^{*} \quad v\left( \det(ae_1,\ldots, ae_n)
  \right)=v\left( \det(e_1,\ldots, e_n) \right) \mod \ n,
\end{equation*}
one can define the \textit{type of a vertex} $x$ in $\bdg$ to be the value
modulo $n$ of the type of one of its representatives. We denote by $\tau(x)$ the
type of $x$.\\
If $L'$ is a second lattice, we can choose our basis $e_1, \ldots, e_n$ for $L$
in such a way that $L'$ admits a basis of the form $a_1e_1, \ldots, a_n e_n$ for
some $a_i \in \bK^*$. The scalars $a_i$ can be taken to be powers of $\pi$. The
incidence relation defined above implies that if the classes of $L$ and $L'$ are
linked by an edge in $\bdg$, then they have different types.
\begin{Rq}\label{Rq:Types}
  Remark that if
  $L=\oplus_{i} \RInt e_i$ and
  \begin{equation*}
    L'=\RInt \pi e_1 \oplus \cdots \oplus \RInt
    \pi e_j \oplus e_{j+1} \oplus \cdots \oplus e_{n},
  \end{equation*}
  then $\tau\big([L\prime]\big)=\tau\big([L\prime]\big)+j\mod \ n$.
\end{Rq}
The action of $SL_n(\bK)$ on $\bdg$ preserves the determinant and is transitive
on the pairs of vertices of the same type. So there are exactly $n$ orbits under
the action of $SL_n(\bK)$ (see \cref{fig:orbitsbuilding} and \cref{fig:ApptPSL3}
for examples).

\paragraph{Apartments}
If $\mbfe$ is a basis of $\bK^n$ then the sub-complex $\calA$ induced by the
set of vertices $\left\{ \oplus^n_{i=1} \RInt \pi^{k_i}e_i \ | \ k_i \in \bZ
\right\}$ is isometric to a $(n-1)$-dimensional Euclidean space tiled by regular
$(n-1)$-simplices. We call such sub-complexes \textit{apartments}. For example an
apartment in the building of $PSL_2(\bQ_2)$ is isometric to $\bR^1$ tiled with
segments of length $1$ (see \cref{fig:ApptArbre}), whereas for $PSL_3(\bQ_2)$ the
apartment are isometric to $\bR^2$ and tiled with triangles (see
\cref{fig:ApptPSL3}). 

For any two points in $\hat{\bdg}$ there exists an apartment containing them. If
$x,y\in \hat{\bdg}$ let $\calA$ be an apartment containing $x$ and $y$ and
define $d_{\hat{\bdg}}(x,y)$ to be equal to the euclidean distance
$d_{\calA}(x,y)$. This definition does not depend on the choice of apartment
$\calA$ and thus endows $\hat{\bdg}$ with a well defined distance. Moreover, this distance
verifies the \emph{negative curvature inequality:} for all $x,y,z\in
\hat{\bdg}$ and $t\in [0,1]$
\begin{equation}\label{eq:NegativeCurvatureInequality}
  d^2_{\hat{\bdg}}(z,tx+(1-t)y)\leq t d^2_{\hat{\bdg}}(z,x) +
  (1-t) d^2_{\hat{\bdg}}(z,y)-t(1-t) d^2_{\hat{\bdg}}(x,y).
\end{equation}
Denote by $d_{\bdg}$ the distance on the $1$-skeleton $\bdg$ assigning length
$1$ to an edge. Then $d_{\bdg}(x,y)$ is greater than $d_{\hat{\bdg}}(x,y)$ for
all vertices $x$ and $y$ in $\bdg$.  

\paragraph{Contractibility} Using the above inequality one can show that the
building is contractible (see \cite{AB} for more details). We can actually show that
\emph{convex sets} 
in $\hat{\bdg}$ are themselves contractible.

\begin{Claim}\label{Claim:ContractibilityofBdg}
  Let $r>0$. Any convex set in $\hat{\bdg}$ is contractible.
\end{Claim}
\begin{proof}
  Let $r>0$ and $\calC$ a convex set in $\hat{\bdg}$ and endow it with the distance
  induced by $d_{\hat{\bdg}}$. Take $x_0\in \calC$ and define,
  \begin{equation*}
    \mathcal{H} : \begin{cases}
      [0,1] \times \calC & \rightarrow \calC,\\
      (t,x) & \mapsto tx+(1-t)x_0.
    \end{cases}
  \end{equation*}
  Since $\calC$ is convex, the map $\mathcal{H}$ is well-defined.
  Moreover $\mathcal{H}(0,\cdot)=id_{\calC}$ and $\mathcal{H}(1,x)=x_0$
  for all $x$ in $\calC$. Let us show that $\mathcal{H}$ is continuous. Take $x,x^\prime \in
  \calC$ and $t,t^\prime \in [0,1]$ and let $z=t^\prime x^\prime +
  (1-t^\prime)x_0$. By \cref{eq:NegativeCurvatureInequality}
  \begin{equation}\label{eq:NCIProof}
    d^2_{\hat{\bdg}}(z,tx+(1-t)x_0)\leq t d^2_{\hat{\bdg}}(z,x) +
    (1-t) d^2_{\hat{\bdg}}(z,x_0)-t(1-t) d^2_{\hat{\bdg}}(x,x_0).
  \end{equation}
  But if $\calA$ is a an apartment containing $z$ and $x_0$, then by property of
  the Euclidean distance $d_{\calA}$
  \begin{equation*}
    d_{\hat{\bdg}}(z,x_0)=d_{\calA}(t^\prime x^\prime +(1-t^\prime)x_0,x_0)
    =t^\prime d_{\calA}(x^\prime,x_0)=t^\prime d_{\hat{\bdg}}(x^\prime,x_0),
  \end{equation*}
  which tends to $td_{\hat{\bdg}}(x,x_0)$ as $(t^\prime,x^\prime)$ tends to
  $(t,x)$. Similarly 
  \begin{align*}
    d_{\hat{\bdg}}(z,x) 
    \leq d_{\hat{\bdg}}(z,x^\prime)
      + d_{\hat{\bdg}}(x^\prime,x)
    &= d_{\hat{\bdg}}(t^\prime x^\prime +(1-t^\prime)x_0,x^\prime) + d_{\hat{\bdg}}(x^\prime,x),\\
    &=(1-t^\prime)d_{\hat{\bdg}}(x^\prime,x_0)+ d_{\hat{\bdg}}(x^\prime,x),
  \end{align*}
  which converges to $(1-t)d_{\hat{\bdg}}(x,x_0)+ d_{\hat{\bdg}}(x^\prime,x)$ as
  $(t^\prime, x^\prime)$ tends to $(t,x)$. Thus the right term of
  \cref{eq:NCIProof} converges to $0$ as $(t^\prime, x^\prime)$ tends to $(t,x)$.
  Hence the continuity of $\mathcal{H}$ and the contractibility of
  $\calC$.
\end{proof}

\subsection{Large scale simple connectedness}
\label{ksc}
For a graph $\calG$ and $k \in \bN$, we define a $2$-complex, noted $P_k(\calG)$,
such that:
\begin{itemize}
\item Its $1$-skeleton is given by $\calG$ ;
\item Its $2$-skeleton is composed of $m$-gons (for $m\in [0,k]$) defined
  by the simple loops of length $m$ in $\calG$ (up to cyclic permutations).
\end{itemize}

\begin{Def}
  We say that $\calG$ is \textbf{k-simply connected} or simply connected at scale
  $k$ if $P_k(\calG)$ is simply connected.
\end{Def}
\begin{Ex}
  Let $G$ be a finitely generated group and $T$ a finite symmetric generating
  set. The Cayley graph $(G,T)$ is simply connected at scale $k$ if and only if $G$
  has a presentation $\langle{T,\mathcal{R}}\rangle$ with relations of length at
  most $k$.
\end{Ex}
\begin{Ex}
 Let $n\geq 2$. The Bruhat-Tits building of $PSL_n(\bK)$ is simply connected at scale $3$.
\end{Ex}
\begin{Rq}
  If $k \leq k^\prime$, then every $k$-simply connected graph is
  $k^\prime$-simply connected.
\end{Rq}

The following proposition allows us to restrict the study of the
LG-rigidity of a graph $\calG$ to some smaller class of graphs.

\begin{Prop}[de la Salle, Tessera, {\cite[Proposition 1.5]{DlST2016}}]\label{Prop:15DlST2}
 Let $k\in \bN$ and $\calG$ be a $k$-simply connected graph, with cocompact
 isometry group. Then $\calG$ is LG-rigid if and only if there exists $R$
 such that every $k$-simply connected graph which is $R$-locally $\calG$ is
 isometric to $\calG$. 
\end{Prop}

To apply this result to our proof we will need to show that the studied graph
$\myg$ is simply connected. The following proposition shows that being simply
connected is invariant under quasi-isometry.

\begin{Prop}[de la Salle, Tessera, {\cite[Theorem 2.2]{DlST2016}}]\label{Prop:Th22DlST2}
  Let $k \in \bN^*$ and let $\calG$ be a $k$-simply connected graph. If $\calH$
  is quasi-isometric to $\calG$, then there exists $k^\prime \in \bN^*$ such that
  $\calH$ is simply connected at scale $k^\prime$. 
\end{Prop}

\label{LGRtopo}
Before moving to the next section, let us mention a consequence of
\cref{Prop:15DlST2}. Indeed, this result allows us to look at the LG-rigidity
notion with a topological point of view. 
Let’s denote $\mathfrak{G}_k$ the set of 
isometry classes of locally finite $k$-simply connected graphs. We can
define a distance on this set by: 
\begin{equation*}
  d_{\mathfrak{G}_k}(X,Y):=\inf \left\{2^{-R} \ : \
    X \ \text{and} \ Y \ \text{are} \ R\text{-close} \right\},
\end{equation*}
which endows $\mathfrak{G}_k$ with a topology. \cref{Prop:15DlST2} implies that a
graph is LG-rigid if and only if its isometry class in $\mathfrak{G}_k$ is
isolated for this topology.

\subsection{Extension of isometries} \label{subsec:Extofisom}
In order to build the “hybrid” graph mentionned above, we will need some result
to extend globally our local definition of edges. We recall here the result
of de la Salle and Tessera \cite[Lemma 4.1]{DlST2015} that will serve our purpose.
% ==================================================
% ------------------ Newpage forced
% ==================================================
\newpage

\begin{Prop}[de la Salle, Tessera] Let $\calG$ be a graph with cocompact 
  isometry group. Given some $r_1 \geq 0$, there exists $r_2>0$ such that:
  for every $g \in \calG$, the restriction to $B_{\calG}(g,r_1)$ of an isometry $f$ :
  $B_{\calG}(g,r_2) \rightarrow \calG$ coincides with the restriction of an element of
  $\Is(\calG)$.
  \label{Lmm:41DlST}
\end{Prop}

It is however not necessarily true that $f$ coincides on the whole $B(g,r_2)$
with an isometry of $\calG$. Indeed, truncating the entire graph to some ball might
allow some kind of flexibility near the boundary of the ball (see
Example \ref{Ex:ContreExLmm41DlST} and \cref{fig:IllustrationLmm41DlST}). Hence, in
order to coincide with a global isometry we need to restrict the local isometry
$f$ to a smaller ball which does not contain the flexible area.

\begin{Ex}\label{Ex:ContreExLmm41DlST}
 Let $\calG$ be the Cayley graph of $\bZ^2$ endowed with its usual generating part.
 We consider in
 \cref{fig:IllustrationLmm41DlST} an isometry $f$ defined on
 $B\big((0,0),1 \big)$ such that $f$ fixes $(0,0)$, $(-1,0)$ and
 $(0,-1)$ (represented by the blue vertices) and exchange $(1,0)$ with $(0,1)$
 (the orange and brown vertices).
 Then $f$ is an isometry from $B\big((0,0),1 
 \big)$ to $B\big((0,0),1 \big)$, but can not coincide with a global isometry
 of $\calG$ on that ball. Indeed, if such a global isometry existed, then it
 should send the vertex $(-1,1)$ (represented by the light brown vertex on the
 left part of the figure) at distance $1$ from both $f(-1,0)=(-1,0)$ and
 $f(0,1)=(1,0)$. Which is impossible since the only point at distance $1$ from
 $(1,0)$ and $(-1,0)$ is $(0,0)$ and it is already the image of $(0,0)$. 
\end{Ex}
\begin{figure}[htbp]
  \centering
  \begin{tikzpicture}[scale=0.9]
  \foreach \x in {-1,0,1}{
    \draw[dashed, black!60] (\x,1.6) -- (\x,-1.6);
    \draw[dashed, black!60] (1.6,\x) -- (-1.6,\x);
  };
  \node[draw, MarCB, circle, fill=MarCB] at (-1,1) {};
  \foreach \x in {(1,0), (0,1), (-1,0), (0,-1)}{
    \draw \x -- (0,0);
  };
  \node[bob1] at (-1,0) {};
  \node[bob1] at (0,-1) {};
  \node[bob1] at (0,0) {};
  \node[draw, circle, fill=MyOrange] at (0,1) {};
  \node[draw, circle, fill=MFCB] at (1,0) {};
  % Arrow of f
  \draw[->,>=latex] (2.5,0) -- (3.5,0) node[midway, above] {$f$};
  %
  % Right part
  \begin{scope}[shift={(6,0)}]
    \foreach \x in {-1,0,1}{
      \draw[dashed, black!60] (\x,1.6) -- (\x,-1.6);
      \draw[dashed, black!60] (1.6,\x) -- (-1.6,\x);
    };
    \foreach \x in {(1,0), (0,1), (-1,0), (0,-1)}{
      \draw \x -- (0,0);
    };
    \node[bob1] at (-1,0) {};
    \node[bob1] at (0,-1) {};
    \node[bob1] at (0,0) {};
    \node[draw, circle, fill=MFCB] at (0,1) {};
    \node[draw, circle, fill=MyOrange] at (1,0) {};
  \end{scope}
\end{tikzpicture}
  \caption{Local isometry that can not coincide with a global isometry on its
    entire domain of definition}
  \label{fig:IllustrationLmm41DlST}
\end{figure}
\subsection{\texorpdfstring{Preliminary results on $\myg$}{Preliminary results
    on the studied graph}}

\begin{Lmm}
 \label{Lmm:PSLinclus} If $\myg$ verifies the hypothesis of
 \cref{Th:BigRigidity}, then $PSL_n(\bK)$ is 
 included in $\rho(\Is(\myg))$. Moreover, if $q(\myg)$ contains a vertex of a
 certain type $i$, then $q(\myg)$ contains all the vertices of type $i$.
\end{Lmm}
\begin{proof}
  Since $\rho\left(\Is(\myg)\right)$ is of finite index in the isometry group of
  the building $\bdg$, the same
  goes for its
  normal core $\cap_{g \in \Is(\bdg)}
  g\rho\left(\Is(\myg)\right)g^{-1}$. Then, by 
  simplicity of $PSL_n(\bK)$, the normal core of $\rho\left(\Is(\myg)\right)$
  contains $PSL_n(\bK)$. Hence the result.

  Then, the second part of the lemma follows from the equivariance of $q$ and
  the transitivity of $PSL_n(\bK)$ on vertices of the same type.
\end{proof}
Without loss of generality, we can assume that $im(q)$ contains type $0$
vertices, that is to say $\tau^{-1}(0) \subset im(q)$. Moreover, using
\cref{Prop:Th22DlST2} we obtain that $\myg$ is simply connected at some scale
$k>0$.  
\vspace{0.2cm}
\begin{center}
  \adforn{11}
\end{center}
\vspace{0.2cm}
The aim of the next two sections is to prove \cref{Th:BigRigidity} for $n \geq
4$. For the sake of clarity we recapitulate here the needed assumptions for the
proof. 
\begin{HypBox}
\textbf{Hypothesis (H)}
\begin{enumerate}
\item Let $\myg$ be a $k$-simply-connected transitive graph;
\item Let $\mY$ be a graph $R$-locally $\myg$ and $k$-simply connected;
\item Let $n\geq 4$ and $\bK$ a non-Archimedean local skew field of characteristic
  zero. Denote by $\bdg$ the Bruhat-Tits building
  of $PSL_n(\bK)$;
\item Let $\rho$~: $\Is(\myg) \rightarrow \Is(\bdg)$ be an injective homomorphism
  and $q$~: $\myg \rightarrow \bdg$ an $\Is(\myg)$-equivariant injective
  quasi-isometry;
\item Assume that $\rho(\Is(\myg))$ is of finite
  index in $\Is(\bdg)$ and that $q(\myg)$ contains $\tau^{-1}(0)$. 
\end{enumerate}
\end{HypBox}

% Section 2
\section{Tracking vertices through their prints}
This section is dedicated to the definition of a graph locally the same as $\bdg$ which we will call
$\CalY$. Before moving to the detailed definition let us explain the idea of the
construction. Recall that the vertices of $\bdg$ are partitioned into different
types (see \cref{BTBuilding}) denoted by integers in $\{0, \ldots, n-1\}$. By
\cref{Lmm:PSLinclus}, if $q(\myg)$ contains a vertex of a certain type then
it contains all the vertices of that type. Denote by $\Tau$ the set of types that
are not contained in $q(\myg)$, namely $\Tau=\{0, \ldots, n-1\}
\backslash \tau(q(\myg))$. We have the following partition
\begin{equation}\label{eq:Partitonbdg}
  \bdg = q(\myg) \sqcup \left(\sqcup_{i \in \Tau} \,
    \tau^{-1}(i)  \right). 
\end{equation}
\begin{Ex}
  Take $\bK=\bQ_2$ and assume that $im(q)$ is composed only of type zero vertices.\\
  When $n=2$ we have $T=\{1\}$ and the
  building is represented in \cref{fig:orbitsbuilding}. The partition
  in \cref{eq:Partitonbdg} corresponds to the partition of vertices in two
  different colors.\\ 
  When $n=3$, we get $T=\{1,2\}$. An apartment of $\bdg$ is represented in
  \cref{fig:ApptPSL3} and the partition of this part of $\bdg$ corresponds to
  the partition in three different colors.
\end{Ex}
\begin{Ex}\label{Ex:partitiondeX4}
  Let $n=4$ and $\bK=\bQ_2$ and assume that $im(q)$ contains type zero and type $2$
  vertices. Then $T=\{1,3\}$. We will not try to represent $\bdg$ or an apartment but
  recall that it is tiled by tetrahedrons. The partition is illustrated on a
  tetrahedron in \cref{fig:tetrahedre}, where $im(q)$ corresponds to the two blue vertices.
\end{Ex}
\begin{figure}[htbp]
  \centering
  \begin{tikzpicture}
        \node[cercle, fill=MarCB] (x0) at (50:1.2) {};
        \node[cercle, fill=Leturquoise] (x1) at (0,0) {};
        \node[cercle, fill=MyOrange] (x2) at (20:2) {};
        \node[cercle, fill=Leturquoise] (x3) at (-35:1.2) {};
        \draw[dashed] (x0) -- (x3);
        \draw (x2) -- (x0) -- (x1) -- (x3) -- (x2);
        \draw (x1) -- (x2);
        \node[anchor=west] (type3) at (25:2.5) {\small Type $3$};
        \node[anchor=south] (type1) at (100:2) {\small Type $1$};
        \node[anchor=east] (type0) at (-180:1.8) {\small $im(q)$ (Type $0$ and $2$)};
        \draw[->,>=latex, black!40] (type3.south) to[bend left] (x2.south east);
        \draw[->,>=latex, black!40] (type1.east) to[bend left] (x0.north);
        \draw[->,>=latex, black!40] (type0.east) to[bend left] (x1.west);
        \draw[->,>=latex, black!40] (type0.south) to[bend right] (x3.south);
  \end{tikzpicture}
  \caption{Partition of a simplex}
  \label{fig:tetrahedre}
\end{figure}
The idea of the construction of $\CalY$ is to take the vertices of $\mY$ and add
to them vertices of the missing types, \textit{ie.} vertices with type in
$\Tau$ (see \cref{fig:SommetsdecalYexemple} for an example). But we want to
build this vertices only with informations encoded in $V(\mY)$.  
That is why we introduce the local characterization of a vertex in the building
(see \cref{sec:localcharacterizationbdg}). Then, using a well chosen set of
isometries from $\mY$ to $\myg$, we transfer this print notion to
$\mY$, each print in $\mY$ corresponding to a vertex of a missing type.

\subsection{Prints in a building}
\label{sec:localcharacterizationbdg}

In this section we show that a vertex in $\bdg$ can be determined by a part of its
$1$-neighbourhood. More precisely, we prove that a vertex in the building is
entirely determined by the vertices in its $1$-neighbourhood having
type zero.

\begin{Def}\label{Def:Print}
 Let $x$ be a vertex of $\bdg$. We define the \textbf{print} of $x$, denoted by
 $\Print(x)$, to be the intersection of the 
 $1$-neighbourhood of $x$ with the vertices of type zero, \textit{viz.}
 $\Print (x):= B_{\bdg}(x,1)\cap \tau^{-1}(0)$.
\end{Def}
\begin{Rq}
  We choose to define print as a set of vertices of type zero because (in order
  to simplify notations and proofs) we assumed from the beginning that
  $\tau^{-1}(0)$ was contained in $im(q)$. But we could have
  taken any other type.
\end{Rq}
\begin{Ex} \Cref{fig:printinX} represents a ball of radius $1$ in two different
  cases. The case when
  $n=2$ and $|\RInt/ \pi \RInt|=2$ (for example when $\bK=\bQ_2$) is represented
  on the left figure. The case 
  when $\bK=\bQ_2$ and $n=3$ is represented on the right figure. In 
  each case, the print of $x$ corresponds to the set of blue vertices.  
\end{Ex}
\begin{figure}[htbp]
  \centering
  \begin{tikzpicture}
    \node[bob2] (x) at (0,0) {};
    \node[bob1] (z) at (0:1) {};
    \node[bob1] (c) at (120:1) {};
    \node[bob1] (d) at (240:1) {};
    \draw (z) -- (x) -- (c);
    \draw (x) -- (d);
    \node[anchor=south] (x1) at (c.north) {$x_1$}; 
    \node[anchor=west] (x3) at (z.east) {$x_3$}; 
    \node[anchor=north] (x2) at (d.south) {$x_2$}; 
    \node[anchor=east] (xx) at (x.west) {$x$};
    \begin{scope}[shift={(0,2)}]
      \node[anchor=south] (Px) at (0,0) {$\Print(x)=\{x_1, x_2, x_3\}$};
    \end{scope}
    \begin{scope}[shift={(0,-2.5)}]
      \node () at (0,0) {$B(x,1)$ for $n=2=p$};
    \end{scope}
    \begin{scope}[shift={(4,0)}]
      \draw[very thin, black!40] (-0.5,-0.5) rectangle (3,1.5);
      \node[bob1, inner sep=3pt] (1) at (0,0) {};
      \node[anchor=west] (l1) at (1.east) {\small Type $0$ vertices}; 
      \node[bob2, inner sep=3pt] (2) at (0,0.5) {};
      \node[anchor=west] (l2) at (2.east) {\small Type $1$ vertices}; 
      \node[draw,circle, fill=MyOrange, inner sep=3pt] (3) at (0,1) {};
      \node[anchor=west] (l3) at (3.east) {\small Type $2$ vertices}; 
    \end{scope}
    \begin{scope}[shift={(10,0)}]
      \node[bob2] (x) at (0,0) {$x$};
      \foreach \y in {1,3,5,7,9,11,13}{
        \draw ({\y *360/14}:1.5)-- (x);
        \draw[fill=Leturquoise, thick] ({\y *360/14}:1.5) circle (0.15);
      }
      \foreach \y in {0,2,...,12}{
        \draw ({\y *360/14}:1.5)-- (x);
        \draw[fill=MyOrange, draw] ({\y *360/14}:1.5) circle (0.15);
      }
      \begin{scope}[shift={(0,-2.5)}]
        \node () at (0,0) {$B(x,1)$ for $p=2$ and $n=3$};
      \end{scope}
    \end{scope}
  \end{tikzpicture}
  \caption{Prints and $1$-neighbourhood of a vertex in $\bdg$}
  \label{fig:printinX}
\end{figure}

The following result proves that a vertex in $\bdg$ is uniquely determined by
its print.

\begin{Prop}\label{Prop:characterisationparprint}
  Let $x_1, x_2 \in \bdg$.
  If $\Print(x_1)=\Print(x_2)$, then $x_1=x_2$.
\end{Prop}
Before showing the above property, let us recall (and prove) a useful
fact concerning the choice of representative of a vertex.
\begin{Claim} \label{Claim:choixdescoordonnees}
  For any vertex in $\bdg$, we can always find a representative $\oplus_i \RInt
  \pi^{k_i}e_i$ of the vertex such that 
  \begin{equation}\label{eq:choixdescoordonees}
    \begin{cases}
      \forall i \in \{1, \ldots, n\} &k_i \geq 0,\\
      \exists i_0 \in \{1, \ldots, n\} &k_{i_0} = 0.
    \end{cases}
  \end{equation}
\end{Claim}
\begin{proof}[Proof of the claim]
  Indeed, let $x\in \bdg$ and let $(l_1, \ldots, l_n)$ be a representative of
  $x$ and let $i_0$ be such that $l_{i_0}=\min_i l_i$, then
  \begin{equation*}
    \left[\oplus^n_{i=1} \RInt {\pi}^{l_i}e_i\right]
    = {\pi}^{- l_{i_0}} \left[\oplus^n_{i=1} \RInt {\pi}^{l_i-l_{i_0}}e_i\right]
    = \left[\oplus^n_{i=1} \RInt {\pi}^{l_i- l_{i_0}}e_i\right].
  \end{equation*}
  Thus $(l_1-l_{i_0}, \ldots, l_1-l_{i_0})$ is a representative of $x$ and
  verifies \cref{eq:choixdescoordonees}.
\end{proof}

Now, let us prove that the print determines the vertex.

\begin{proof}[Proof of {\cref{Prop:characterisationparprint}}]
  Let $x_1, x_2 \in \bdg$ such that $\Print(x_1)=\Print(x_2)$. 
  
  First remark that if $\tau(x_1)=0$ then $\Print(x_1)=\{x_1\}$ which implies that
  $\Print(x_2)=\Print(x_1)=\{x_1\}$. But then $x_2$ has only one
  neighbour of type $0$, which is only possible if $\tau(x_2)=0$. Thus
  $\{x_2\}=\Print(x_2)=\{x_1\}$ and so $x_1=x_2$. 

  Now assume that $\tau(x_1)\neq 0$ and
  take $\Appart$ to be an apartment containing both $x_1$ and $x_2$. Define
  $P:=\Print(x)\cap \Appart$ and let $\mbfe$ be a basis such that
  \begin{equation*}
    \Appart=\left\{ \oplus^n_{i=1} \RInt {\pi}^{k_i}e_i \ | \ k_i \in \bZ \right\}
    \quad \text{and} \quad
    x_1=(0, \ldots,0).
  \end{equation*}
  By \cref{Claim:choixdescoordonnees}, we can choose a representative
  $(k_1,\ldots,k_n)$ of $x_2$ such that $k_i\geq 0$ 
  for all $i$ and there exists $j\in \{1,\ldots, n\}$ such that $k_{j}=0$.
  Now define the sequence $i_1, \ldots, i_n$ of indices such that 
  $k_{i_n} \geq \cdots \geq k_{i_1}=0$ and let
  \begin{equation*}
    l_{i_1}= \cdots = l_{i_{\tau(x)}} =1 \qquad
    l_{i_{n-\tau(x)+1}}= \cdots = l_{i_n}=0.
  \end{equation*}
  Then by \cref{Rq:Types} the vertex $z=(l_1, \ldots,l_n)$ has type $0$.
  Moreover it is at distance $1$ from $x_1$, so $z$ belongs to $P$. But if
  $k_{i_n}>0$, then $d(z,x_2)>1$ thus $z$ can not belong to 
  $\Print(x_2)$. Hence $k_{i_1}\leq 0$, that is to say $k_i=0$ for all $i$ and
  thus $x_2=x_1$.
\end{proof}

This proves that a vertex in $\bdg$ is uniquely determined by its print. Thus,
we can introduce the following definition without ambiguity. 

\begin{Def}
  Let $x$ to be a vertex in $\bdg$. We say that $x$ is the
  \textbf{source} of $\Print(x)$.
\end{Def}

In order to prove \cref{Th:BigRigidity}, we will need to know how prints behave
under the action of $PSL_n(\bK)$. So let $x \in \bdg$ and let $\alpha \in
PSL_n(\bK)$. Since $\alpha$ is an isometry, we get
\begin{equation*}
  \alpha \Big(\Print(x) \Big)= \alpha \Big( B(x,1) \cap \tau^{-1}(0) \Big)
  =  \alpha \big(B(x,1)\big) \cap \alpha\tau^{-1}(0)
  = B\left( \alpha(x),1 \right) \cap \tau^{-1}(0).
\end{equation*}
We deduce the following lemma.
\begin{Lmm} \label{Lmm:printinvparPSL}
  Let $x \in \bdg$. If $\alpha$ belongs to $ PSL_n(\bK)$, then $\alpha \left( \Print(x)
  \right)= \Print \left( \alpha(x)\right)$. 
\end{Lmm}

\subsection{Atlas of local isometries}
To build our graph locally the same as $\bdg$, we need to restrict ourselves to
a particular set of local isometries from $\mY$ to $\myg$. More precisely, if
$y_1$ and $y_2$ are close in $Y$ and $f_1$ (resp. $f_2$) is an isometry
from $B_{\mY}(y_1,R)$ (resp. $B_{\mY}(y_2,R)$) to $\myg$, we
want the transition map $f_2f^{-1}_1$ to coincide with an element in
$\rho^{-1}PSL_n(\bK)$ on a small ball. This is what we formalize here and
schematize in \cref{fig:constructiondelatlas}.

In order to avoid any ambiguity regarding the notion of center of a ball,
let us precise our definition of ball in a graph. What we mean when we talk
of “a ball of radius $R$” is actually a \emph{pointed ball of
  radius $R$} that is to say, a couple $(\mathcal{B},y)$ such that $y$ is a
vertex in $\mY$ and $\mathcal{B}=B_{\mY}(y,R)$. We will abuse notation by denoting
such a pointed ball $B_Y(y,R)$ (instead of $\big(B_Y(y,R), y\big)$). This way,
the center of a ball is always well defined. 

\begin{Def}
  Let $\Atlas$ be a set of isometries from balls of radius $R$ in $\mY$ to
  $\myg$. We say that $\Atlas$ is an \textbf{atlas} of local
  isometries from $\mY$ to $\myg$ if
  the map that associates to each isometry in $\Atlas$ the center of its
  ball of definition is a bijection from $\Atlas$ to $Y$. That is to say, we
  can write 
  \begin{equation*}
    \Atlas:= \left\{ f_y \ : \ B_{\mY}(y,R) \rightarrow \myg \ | \ y \in \mY \right\},
  \end{equation*}
  where the map that associates $f_y$ to $y$ is bijective.\\
  We say that $f_y$ is the \textbf{isometry associated to $y$} in $\Atlas$.
\end{Def}

Let $H_0:=\rho^{-1}PSL_n(\bK) $. Now, we show that we can construct an atlas
of local isometries from $\mY$ to
$\myg$ such that the transition maps between two isometries defined on balls
with neighbouring centers coincide with elements of $H_0$.

We will note a path
between two vertices $v_1$ and $v_2$ as a sequence $(v_1,\ldots, v_l)$ of
adjacent vertices.
% =========================================
% ------------ Newpage forced
%=========================================
\newpage

\begin{Lmm}
  \label{Prop:ConstructionDeLAtlas}
  Let $\ratlas >0$ and let $H_0:= \rho^{-1}PSL_n(\bK)$. For $R$ large enough, if
  $\mY$ is $R$-locally $\myg$, then there exists an atlas $\Atlas$ such that
  for any two neighbours $y$ and $z$ in $\mY$
  \begin{equation}
    \label{eq:transitionmap}
    \exists a \in H_0 \quad
    {f_y f^{-1}_z}_{|B(f_z(z),\ratlas)} = a_{|B(f_z(z),\ratlas)}.
  \end{equation}
\end{Lmm}
Before proving it, let us schematize the framework of this lemma. 
In \cref{fig:constructiondelatlas} we represent two isometries $f_y$ and $f_z$
with $z$ neighbour to $y$. The larger discs correspond to balls of radius $R$
and the smaller ones to balls of radius $\ratlas$.
The map $f_yf^{-1}_z$ restricted to $B(f_z(z),\ratlas)$ takes $f_z(z)$ to
$f_y(z)$ which is a neighbour of $f_y(y)$ and coincides on this ball with an
element in $H_0$.

\begin{figure}[htbp]
  \centering
  \begin{tikzpicture}[scale=0.9]
    \node (Y) at (-2.5, 0) {$\mY$};
    % Gros cercle bleu
    \draw[fill=Leturquoise!70, opacity=0.5] (0,0) circle (1.5);
    % Les cercles marrons
    \begin{scope}[shift={(0.4,0)}]
      \draw[fill=MFCB!70, opacity=0.5] (0,0) circle (1.5) ;
      \draw[fill=MFCB, opacity=0.7] (0,0) circle (0.8);
    \end{scope}
    % Nœud étiqueté y
    \draw[fill=black] (0,0) circle (0.05) node[above](y) {$y$};
    % ---------------- Partie de droite ---------------
    \begin{scope}[shift={(8,0)}]     
      \node (X) at (3, 0) {$\myg$};
      % Gros cercle bleu de droite
      \draw[fill=Leturquoise!70, opacity=0.5] (0,0) circle (1.5);
      \begin{scope}[shift={(0.5,-2)}]
        % Les cercles marrons 
        \draw[fill=MFCB!70, opacity=0.5] (0,0) circle (1.5) ;
        \draw[fill=MFCB, opacity=0.7] (0,0) circle (0.8);
        \node (ptfzm1) at (-80:0.8){};
        \node (ptarrivee) at (180:1.5){};
        \node (dept) at (40:0.8) {}; % Pt de dpt de la fleche courbe
        \draw[fill=black] (0,0) circle (0.05) node[below,anchor=north west](fzz) {\scriptsize $f_z(z)$};
      \end{scope}
      % Le petit cercle turquoise
      \begin{scope}[shift={(0,0.4)}]
        \node (blob) at (-180:0.8){};
        \draw[fill=Leturquoise, opacity=0.7] (0,0) circle (0.8);
        \draw[fill=black] (0,0) circle (0.05) node[above](fyz) {\scriptsize{ $f_y(z)$}};
        \node (arrivefyfz) at (60:0.8) {};
      \end{scope}
      \draw[fill=black] (0,0) circle (0.05) node[below, anchor=north east](fyy) {\scriptsize{$f_y(y)$}};
      \draw[thin] (0,0.4) -- (0,0);
      % La fleche de fyfz{-1}
      \draw[->, >=latex, thick] (dept.center) to[bend right] (45:1.9)%
      to[bend right]  (arrivefyfz.center);
      \node[] (fyfz) at (45:2.3) {$f_yf^{-1}_z$};
    \end{scope}
    % Flèches entre
    \begin{scope}[shift={(0.5,0)}]
      \draw[fill=black] (0,0) circle (0.05) node[below](z) {$z$}; % Nœud z
      \draw[->, >=latex, MFCB] (-30:1.4) to[bend right] (ptarrivee.center); % Grande fleche marron
      \draw[<-, >=latex, MFCB] (-90:0.8) to[bend right] (ptfzm1.center); % Petite fleche marron
      \draw[->, >=latex, Leturquoise] (0:0.70) to[bend left] (blob.center);% La petite fleche bleue
    \end{scope}
    \draw[->, >=latex, Leturquoise] (100:1.5) to[bend left] +(8,0);
    \draw[thin] (0,0) -- (0.5,0);
    % Les labels de flèches
    \node[] (fy) at (4,2.9) {$f_y$};
    \node[] (fz) at (4,-1.7) {$f_z$};
    \node[] (fzm) at (4,-3.6) {$f^{-1}_z$};
    \node[] (fym) at (4,0.4) {$f_y$};
  \end{tikzpicture}
  \caption{Composition of isometries with neighbouring centers}
  \label{fig:constructiondelatlas}
\end{figure}

Let us discuss the idea of the proof. First, for two neighbours $y$ and $z$ we
use \cref{Lmm:41DlST} to prove 
that $f_yf^{-1}_z$ coincides on a small ball with an element $a$ in $\Is(\myg)$.
This isometry corresponds to the “default” of belonging to $H_0$ we want to
correct.
Hence, we consider in our atlas the new isometry defined on $B(z,R)$ by
$af_z$. Finally, we extend this construction along paths in $\mY$ and prove
that the wanted property for $\mathfrak{A}$ does not depend on the choice of
path.

\begin{proof} Let $\ratlas >0$ and let $H_0:= \rho^{-1}PSL_n(\bK)$.
  Now, let $y\in \mY$ and $f_y$ be an isometry from $B(y,R)$ to $\myg$. Let $z$ be a
  neighbour of $y$ in $Y$ and $\tilde{f}_z$ be an isometry from $B(z,R)$ to $\myg$. Then
  the map
  \begin{equation*}
    f_y {\tilde{f}}^{-1}_z \ : \ B_{\myg}\Big({\tilde{f}}_z(z),R-1\Big)
    \rightarrow  B_{\myg}\left( f_y(z),R-1 \right)
  \end{equation*}
  is a well defined local-isometry of $\myg$. By \cref{Lmm:41DlST} if $R$ is large
  enough, there exists $a$ in $\Is(\myg)$ such that $ f_y
  {\tilde{f}}^{-1}_z$ coincides with $a$ on $B_{\myg}\left(\tilde{f}_z(z),\ratlas +k\right)$, where we
  recall that $k$ refers to the scale at which $\mY$ is simply connected. We
  will see below why we need to consider such a radius.
  
  Now let ${f}_z:= a\tilde{f}_z$. By definition we have
  \begin{equation*}
    {f}_z \ : \
    \begin{cases}
      B_{\mY}(z,R) &\rightarrow  B_{\myg}\left(f_y(z), R\right),\\
      z & \mapsto a{\tilde{f}}_z(z)=f_y(z),
    \end{cases}
  \end{equation*}
  thus the transition map $f_y{{f}^{-1}}_z$
  is well defined on $B_{\myg}({f}_z(z),R-1)$. Moreover, by choice of $f_z$ we
  get that $f_y{{f}^{-1}}_z$ restricted to $B(f_y(z),\ratlas+k)$ coincides with
  the identity and thus belongs to $H_0$.
  Extending this construction along paths in $\mY$ we get an atlas
  $\Atlas$ of local isometries from $\mY$ to $\myg$.

  Now if $y \in Y$ and
  $f_y$ is the associated isometry in $\Atlas$, we want to show that (up
  to a multiplication by an element in $PSL_n(\bK)$) this isometry 
  does not depend on the choice of path. So let $y \in Y$ and $(y_0=y, y_1, \ldots, y_l=y)$ be a loop
  of length $l$. Take $f_0$ to be an isometry from $B_Y(y_0,R)$ to $\myg$ and
  using the process detailed above, build a sequence of isometries $f_1, \ldots,
  f_l$ such that $f_i$ is defined on $B_Y(y_i,R)$ and
  \begin{equation*}
    \forall i \in \{1, \ldots, l\} \
    \exists a_i \in H_0 \ | \ 
    {\left( f_{i-1}f^{-1}_{i} \right)}_{|B(f_i(y_i),\ratlas +k)}
    = {a_i}_{|B(f_i(y_i),\ratlas +k)}.
  \end{equation*}
  We have to prove that the restrictions to $B(y_0,\ratlas)$ of $f_0$ and $f_l$
  are equal up to a multiplication by an element 
  in $H_0$. Since $Y$ is simply connected at scale $k$, we only have to
  prove this for loops of length smaller than $k$. Hence, we assume that $l \leq
  k$.

  First, remark that for all $i\in \{0, \ldots, l-1\}$

  \begin{equation*}
    \begin{cases}
      f_{i-1}f^{-1}_i \ : \ B_{\myg}\left(f_i(y_i),\ratlas+k\right)
     &\rightarrow B_{\myg}\left(f_{i-1}(y_i),\ratlas +k\right),\\
      f_{i}f^{-1}_{i+1} \ : \ B_{\myg}\left(f_{i+1}(y_{i+1}),\ratlas +k\right)
     &\rightarrow B_{\myg}\left( f_{i}(y_{i+1}),\ratlas+k \right).
    \end{cases}
  \end{equation*}
  Now since $y_i$ and $y_{i+1}$ are at distance $1$, the ball
  $B_{\myg}\left(f_{i}(y_{i+1}),\ratlas+k-1 \right)$ is included in
  $B_{\myg}(f_i(y_i),\ratlas+k)$. Hence the map $\left( f_{i-1}f^{-1}_i  \right)
  \left( f_{i}f^{-1}_{i+1}  \right)$ is well defined and coincides with
  $a_{i}a_{i+1}$ on $B_{\myg}\left(f_{i+1}(y_{i+1}),\ratlas+k -1 \right)$. By
  induction we get that for all $x$ in $B_{\myg}\left(f_{i+1}(y_{i+1}),\ratlas +k-l
  +1\right)$
  \begin{equation*}
    f_0f^{-1}_l(x)=
    \left(f_{0}f^{-1}_1\right) \cdots   \left(f_{l-1}f^{-1}_l\right) (x) = a_1 \cdots a_l(x). 
  \end{equation*}
  Since $\prod^l_{i=1}a_i$ belongs to $H_0$ and $l$ is smaller than
  $k$, it implies that $f_0$ is equal to $f_l$ on $B_Y(y_0,r_A)$ up to
  multiplication by an element in $H_0$.
\end{proof}

The atlas is defined such that a transition map between two isometries defined on
balls with neighbouring centers belongs to $H_0$. But in fact, this property is
also true when the centers are at a slightly bigger distance. 

\begin{Lmm}
  \label{Lmm:transitionmapdefcalY}
  \label{Rq:chgmtdisom}
  Let $r>0$ and $\Atlas$ be an atlas verifiying the conditions of
  \cref{Prop:ConstructionDeLAtlas} with $\ratlas >3r$. Let $y$ and $z$
  in $\mY$ be at distance less than $2r$ and $f_y$, $f_z$ the associated
  isometries in $\Atlas$. Then 
  \begin{equation}
   \label{eq:transitionmapdefcalY}
    \exists a \in H_0 \qquad
    {\left( f_{y}f^{-1}_{z} \right)}_{|B_{\mY}(z,r)}
    = {a}_{|B_{\mY}(z,r)}.
  \end{equation}
\end{Lmm}
\begin{proof}
  Let $r>0$ and assume $\ratlas >3r$. Let $y,z \in \mY$ be at distance $l \leq
  2 r$ and let $f_y,f_z$ be two elements of $\Atlas$ such that
  \begin{equation*}
    f_y \ : \ B_Y(y,R) \rightarrow \myg \quad
    f_z \ : \ B_Y(z,R) \rightarrow \myg.
  \end{equation*} 
  Take $(y_0=y, y_1, \ldots, y_l=z)$ to be a geodesic between $y$ and $z$, and
  for all $i \in \{0, \ldots, l\}$, let 
  $f_i \in \Atlas$ be the isometry associated to $y_i$. Remark that by
  definition of an atlas, it implies $f_0=f_y$ and $f_l=f_z$ and
  \begin{equation*}
    \forall i \in \{0,\ldots,l-1\} \quad \exists a_i \in H_0 \quad 
    {\left(f_i f^{-1}_{i+1}\right)}_{|B\left( f_{i+1}(y_{i+1}),\ratlas\right)}
    = {a_i}_{|B\left( f_{i+1}(y_{i+1}),\ratlas\right)}.
  \end{equation*}
  Now, if $\ratlas > 3 r$ and $l\leq 2r$, then $B_{\mY}(z,r)$ is contained in
  $B_{\mY}(y,\ratlas)$.
  Hence the composition of transition maps $ \left(f_0 f^{-1}_{1}\right) \cdots
  \left(f_{-1} f^{-1}_{r}\right)$ is well defined on $B_{\mY}\left(f_l (y_l),
    \ratlas -l\right)$ and verifies on that ball
  \begin{equation}
    \label{eq:chgmtcarte1}
    f_0 f^{-1}_l=\left(f_0 f^{-1}_{1}\right) \cdots \left(f_{l-1} f^{-1}_{l}\right)
    = a_0 \cdots a_{l-1}.
  \end{equation}
  Hence the result. 
\end{proof}
\subsection{\texorpdfstring{Prints in $Y$}{Prints in Y}}

Using the atlas built above, we can now transfer this print notion to the graph
$\mY$. Let $\rprint>0$ and assume that $\mY$ is endowed with an atlas of isometries
$\Atlas$ as given by \cref{Prop:ConstructionDeLAtlas} with $\ratlas> 3 \rprint$.
Hence, we have
\begin{equation*}
  R>\ratlas > 3 \rprint> \rprint.
\end{equation*}

\begin{Def}\label{Def:PrintY} 
  Let $\PY$ be a set of vertices in $\mY$. 
  We say that $\PY$ is a \textbf{print} if there exists $y$ in $\mY$ and $f
  \in \Atlas$ an isometry from $B_{\mY}(y,R)$ to $\myg$ such that
  \begin{itemize}
  \item The set $\PY$ is contained in $B_{\mY}(y,\rprint)$;
  \item There exists $x \in \bdg \backslash im(q)$ such that $\Print(x)=qf(P)$.
  \end{itemize}
\end{Def}
\begin{Rq}
  Note that in the definition above we ask that $x$ does \emph{not} belong to
  $im(q)$. The definition would also make sense if $x$ belonged to $im(q)$ but
  the purpose of these prints is to reconstruct the ”missing” vertices, namely
  vertices that are not in the image of $q$. Thus to simplify formalism in the
  next pages, we chose to restrict now the definition to prints of vertices in
  $\bdg\backslash im(q)$.   
\end{Rq}
\begin{Ex} If $n=3$ and $p=2$ there are exactly $3$ types of vertices, each
  represented in \cref{fig:PrintinY} by a different color. The $1$-neighbourhood
  of a vertex $x$ in $\bdg$ is then composed of fourteen vertices, represented on
  the right side of the aforementioned figure (where $x$ is the brown vertex at
  the center). 
  If $x\in \bdg\backslash im(q)$ then seven of these
  fourteen vertices are in $im(q)$ (the blue vertices). 
  On the left side of the figure is
  represented $P$ (the black dots) inside $B(y,\rprint)$ (the darker disc).
  The set $qf(P)$ is exactly the set of blue vertices. Hence $P$ is a print.
\end{Ex}
\begin{figure}[htbp]
  \centering
  \begin{tikzpicture}
    % ----------------------------------------- Y
    \begin{scope}[shift={(0,0)}]
      \node (Y0) at (-0.75,0.25) {$\mY$};
      \draw[fill=TurquoiseClair] (1.75,0.25) circle (1.5);
      \draw[fill=Leturquoise!70] (1.75,0.25) circle (1); 
      \draw[fill=black] (1.75,0.25) circle (0.02);
      \draw[fill=black] (2.05,0.25) circle (0.02);
      \draw[fill=black] (2.0,0.75) circle (0.02);
      \draw[fill=black] (1.23,-0.5) circle (0.02);
      \draw[fill=black] (1.75,-0.73) circle (0.02);
      \draw[fill=black] (1.45,0.05) circle (0.02);
      \draw[fill=black] (1.10,0.45) circle (0.02);
      \draw[->, >=latex] (3.75,0.25) -- (4.75,0.25) node[midway, above] {$qf$};
    \end{scope}
    % ----------------------------------------- X
    \begin{scope}[shift={(5,0)}]
      \node (bdg) at (5,0.25) {$\bdg$};
      \draw[fill=MarCB!20] (2,0.25) ellipse [x radius=1.75,y radius=1.6];
      \draw[fill=MarCB] (2,0.25) ellipse [x radius=1.5,y radius=1.25];
      \begin{scope}[shift={(1.85,0.3)}]
        \foreach \y in {1,3,5,7,9,11,13}{
          \draw ({\y *360/14}:0.7)-- (0,0);
          \draw[fill=Leturquoise, thick] ({\y *360/14}:0.7) circle (0.05);
        }
        \foreach \y in {0,2,...,12}{
          \draw ({\y *360/14}:0.7)-- (0,0);
          \draw[fill=MyOrange, draw=MyOrange] ({\y *360/14}:0.7) circle (0.05);
        }
      \end{scope}
      \draw[draw=MFCB, fill=MFCB] (1.85,0.3) circle (0.05);
    \end{scope}
    % % ----------------------- Caption
    \begin{scope}[shift={(-0.5,-2.9)}]
      \draw[very thin, black!20] (-0.5,0.35) rectangle (10.85,-1.35);
      \draw[fill=TurquoiseClair] (0,0) circle (0.1) node [right] {$B(y,R)$};
      \draw[fill=Leturquoise] (0,-0.5) circle (0.1) node [right] {$B(y, \rprint)$};
      \draw[fill=black] (0,-1.0) circle (0.02) node [right] {Elements of $P$};
      \begin{scope}[shift={(5,0)}]
        \draw[fill=MarCB!20] (-0.7,0) circle (0.1) node [right] {$qf\big(B(y,R)\big)$};
        \draw[fill=MarCB] (-0.7,-0.5) circle (0.1) node [right] {$qf\big(B(y,\rprint)\big)$};
        \draw[fill=Leturquoise, thick] (-0.7,-1) circle (0.1) node [right] {$qf(P)=\Print(x)$};
        \draw[fill=MFCB, draw=MFCB] (3,0) circle (0.1) node [right] {$x$};
        \draw[fill=MyOrange, draw=MyOrange] (3,-0.5) circle (0.1) node [right] {$B(x,1) \backslash qf(P)$};
      \end{scope}
    \end{scope}
  \end{tikzpicture}
  \caption{Definition of a print in $\CalY$}
  \label{fig:PrintinY}
\end{figure}

For now, let’s say that $\PY$ verifying the definition above is a print
\textit{associated to $y$ and $f$}. We are going to show that this 
definition depends neither on $y$ nor $f$.

\begin{Lmm}
  \label{Lmm:Printinvparlocisometrie}
  Let  $y_1, y_2 \in \mY$ and $f_1, f_2$ be the associated isometries in
  $\Atlas$. Let $\PY$ be a print associated to $y_1$ and $f_1$.\\
  If
  $\PY \subset B(y_2,\rprint)$ then $\PY$ is a print associated to
  $y_2$ and~$f_2$.
\end{Lmm}
\begin{proof} First, remark that since $P \subset  B(y_2,\rprint)\cap
  B(y_1,\rprint) $, then taking any $y$ in $P$ we get
  \begin{equation*}
    d_{\mY}(y_1,y_2) \leq d_{\mY}(y_1,y) + d_{\mY}(y,y_2) \leq 2 \rprint.
  \end{equation*}
  Applying \cref{Lmm:transitionmapdefcalY} with $r=\rprint$, we get that there
  exists $a \in H_0$ such that
  \begin{equation*}
    {\left(f_{1}f^{-1}_{2} \right)}_{|B_{\myg}(f_2(y_2),\rprint)}
    = {a}_{|B_{\myg}(f_1(y_2),\rprint)}.
  \end{equation*}
  Now let $x \in \bdg$ be such that $\Print(x)=qf_1(P)$. Using
  the equivariance of $q$ and \cref{Lmm:printinvparPSL}, we get
  \begin{equation*}
    qf_2(P)=\rho(a)^{-1}qf_1(P) = \rho(a)^{-1}\Print(x) = \Print\left(\rho(a)^{-1}(x)\right).
  \end{equation*}
  Hence $P$ is a print associated to $y_2$ and $f_2$.
\end{proof}
  
This last lemma proves that being a print does not depend on the
choice of local isometry. 
\begin{Rq}
  In the above proof $\rho(a)^{-1}(x)$ has same type as $x$ since $\rho(a)$
  is type preserving. Thus, once we have taken our atlas in $PSL_n(\bK)$, the
  type of the source of $qf(P)$ does not depend on the choice of local isometry $f$.
\end{Rq}
\subsection{\texorpdfstring{Definition of $\calY$: a building’s replica}{A
    building’s replica}} \label{Sec:buildingY}
The following property defines the graph $\calY$ we will demonstrate to be
locally the same as $\bdg$.

\begin{Prop}
  \label{Prop:DefdeYX} Let $\rprint>0$ and $\Atlas$ be the atlas given by
  \cref{Prop:ConstructionDeLAtlas} for $\ratlas > 3\rprint$. If $R$ is large
  enough, then the following
  graph is well defined. \\
  Let $\CalY$ be the graph whose vertices are given by
  \begin{equation*}
    V(\CalY):=V(Y)\sqcup
    \left\{ P \ : \ \exists x \in \bdg\backslash im(q), \ \Print(x)=P \right\},
  \end{equation*}
  and edges are given by:
  \begin{itemize}
  \item If $y_1, y_2 \in V(\CalY)$, then $(y_1,y_2)$ is an edge if there exists
    $z$ in $\mY$ and $f \in \Atlas$ defined on $B_{\mY}(z,R)$ such that $y_1,y_2
    \in B(z,\rprint)$ and  $d_{\bdg}(qf(y_1), qf(y_2))=1$.
  \item If $y \in V(\CalY)$ and $P$ is a print, then $\big(y,P \big)$ is an edge
    if there exists $z$ in $\mY$ and $f \in \Atlas$ defined
    on $B_{\mY}(z,R)$ cointaining $y$ and $P$ and such that 
    $qf(y)$ is at distance $1$ from the source of $qf(P)$.
  \item If $P_1$ and $P_2$ are two prints, then
    $\big(P_1, P_2\big)$ is an edge if there exists $z$ in $\mY$ and
    $f \in \Atlas$ defined on $B_{\mY}(z,R)$ such that $P_1, P_2 \subset
    B_{\mY}(z,\rprint) $ and such that the source of $qf(P_1)$ is at distance
    $1$ from the source of $qf(P_2)$. 
  \end{itemize}
\end{Prop}
Before looking at the proof of this property, let us sketch some part of this
graph.
\begin{Ex} \label{Ex:sommetsCalY}
  If $n=4$ then $\bdg$ is composed of vertices of type $0$, $1$, $2$ and $3$.
  Assume that
  $q(\myg)$ is composed of vertices of type $0$ and $2$, then $T=\{1,3\}$ and we saw the
  corresponding partition of $\bdg$ in \cref{Ex:partitiondeX4} and \cref{fig:tetrahedre}. The
  appearance of the corresponding $V(\CalY)$ is represented in
  \cref{fig:SommetsdecalYexemple}.
\end{Ex}
\begin{figure}[htbp]
  \centering
  \begin{tikzpicture}
    \fill[MFCB] (0,0.5) -- (-0.5,0) -- (3.5,0) -- (4,0.5) -- cycle;
    \draw[pattern={Lines[angle=65, distance=2pt, line width=0.85pt]},%
          pattern color=MyOrange!80] (0,0.5) -- (-0.5,0) -- (3.5,0) -- (4,0.5) -- cycle;
    \draw (0,0.5) -- (-0.5,0) -- (3.5,0) -- (4,0.5);
    \node[anchor=west] (Y0) at (5,0.25) {Prints}; 
    \begin{scope}[shift={(0,1)}]
      \fill[Leturquoise] (0,0.5) -- (-0.5,0) -- (3.5,0) -- (4,0.5) -- cycle;
      \draw (0,0.5) -- (-0.5,0) -- (3.5,0) -- (4,0.5);
      \node[anchor=west]  (Y0) at (5,0.25) {$V(Y)$};
    \end{scope}
  \end{tikzpicture}
  \caption{Schematic view of $V(\CalY)$ in the case of \cref{Ex:sommetsCalY}}
  \label{fig:SommetsdecalYexemple}
\end{figure}

\begin{proof} Let $\CalY$ be as in \cref{Prop:DefdeYX} and let us show that the definition of
  the edges does not depend on the choice of $f$ in the atlas.

  First, let $y_1, y_2 \in \mY$ and $y,z \in \mY$ such that $y_1$ and $ y_2 $
  belong to $B(y,\rprint) \cap B(z,\rprint)$. Then, take two local maps $f_y,
  f_z$ in $\Atlas$ associated to $y$ and $z$ respectively. Then $d(y,z) \leq 2
  \rprint$ and by \cref{Rq:chgmtdisom} there exists $a \in \Is(\myg)$ verifying
  \cref{eq:transitionmapdefcalY}. Hence, by $\Is(\myg)$-equivariance of $q$ we get
  \begin{align*}
    d_{\bdg} \Big(qf_z(y_1), qf_z(y_2) \Big) 
    &= d_{\bdg} \Big( \rho(a) qf_z(y_1), \rho(a) qf_z(y_2) \Big)\\
    &= d_{\bdg} \Big(q \left( af_z(y_1) \right), q \left( af_z(y_2) \right) \Big)
    = d_{\bdg} \Big(qf_y(y_1), qf_y(y_2) \Big).
  \end{align*}
  Thus $ d_{\bdg} \Big(qf_z(y_1), qf_z(y_2) \Big) =1$ if and only if $d_{\bdg}
  \Big(qf_y(y_1), qf_y(y_2) \Big)=1$ and the definition of edges between two
  vertices of $Y$ does not depend on the choice of local isometry.

  Now take $y \in \mY$ and let $P\subset Y$ be a print. Let
  $z$ and $z^\prime$ such that $y$ and $P$ are contained in
  $B(z,\rprint)\cap B(z^\prime,\rprint)$ and take $f$ (resp. $f^\prime$) in
  $\Atlas$ defined on $B(z,R)$ (resp. $B(z^\prime,R)$). Then $d(z,z^\prime) \leq 2
  \rprint$ and by \cref{Rq:chgmtdisom} there exists $a \in \Is(\myg)$ verifying
  \cref{eq:transitionmapdefcalY}. Hence, 
  \begin{align*}
    d_{\bdg} \Big(qf(y), x \Big) 
    &= d_{\bdg} \Big( \rho(a) qf(y), \rho(a)(x) \Big)\\
    &= d_{\bdg} \Big(q \left( af(y) \right), \rho(a)(x) \Big)
    = d_{\bdg} \Big(qf^\prime (y),  \rho(a)(x) \Big).
  \end{align*}
  If $x$ is the source of $qf(P)$ then, by \cref{Lmm:printinvparPSL} we get
  \begin{equation*}
    \Print(\rho(a)(x)) = \rho(a) \left( \Print(x) \right) = \rho(a) qf(P) = qf^\prime (P). 
  \end{equation*}
  Thus, the existence of en edge between $y$ and $P$ in $\CalY$ does not
  depend of the choice of map in $\Atlas$.

  Finally, take $P_1, P_2 \subset Y$ two prints and let
  $z, z^\prime$ in $\mY$ and $f \in \Atlas$ (resp. $f^\prime$) defined on
  $B_{\mY}(z,R)$ (resp. $B(z^\prime,R)$) such that $P_1, P_2 \subset B_{\mY}(z,\rprint)
  \cap B_{\mY}(z^\prime,\rprint)$. Again $d(z,z^\prime) \leq 2
  \rprint$ and by \cref{Rq:chgmtdisom} there exists $a \in \Is(\myg)$ verifying
  \cref{eq:transitionmapdefcalY}. Hence if $x_1$ is the source of $qf(P_1)$ and $x_2$
  the source of $qf(P_2)$, then $d(x_1,x_2)=1$ if and only if $d(\rho(a)(x_1),
  \rho(a)(x_2))=1$. Moreover, by \cref{Lmm:printinvparPSL}
  \begin{equation*}
    \forall i = 1,2 \quad \Print(\rho(a)(x_i))
    = \rho(a) \left( \Print(x_1) \right) = \rho(a) qf(P_i) = qf^\prime (P_i). 
  \end{equation*}
  Hence the existence of en edge between $P_1$ and $P_2$ in $\CalY$ does not
  depend of the choice of map in $\Atlas$.
\end{proof}

% % Section 3
\section{From one graph to the other} \label{Sec:transferisom}
In this section we prove the isometry between the graph $\CalY$
built and the Bruhat-Tits building and show that it induces an isometry between $\myg$
and $\mY$. 
\subsection{Isometry with the building}
We can now prove that $\calY$ is isometric to the Bruhat-Tits building. Recall that
$\ratlas$ is the radius used to define our atlas $\Atlas$ (see
\cref{Prop:ConstructionDeLAtlas}) and $\rprint$ is the radius used to define
prints in $\CalY$ (see \cref{Def:PrintY}). These constants verify $R > \ratlas
>3\rprint> \rprint$.

\begin{Lmm}\label{Lmm:calYRlocX} Let $\Rbdg>0$.
  If $\rprint$ (and hence $R$) is large enough, then $\CalY$ is ${\Rbdg}$-locally
  $\bdg$. 
\end{Lmm}

To prove this lemma, we define explicitely the local isometries on balls of
radius $\Rbdg$ and prove
that these maps are well defined injections. Then, we compute the minimal value of
$\rprint$ necessary for these applications to be surjective on balls of radius
$\Rbdg$. 
We conclude by showing that these
maps preserve the distance.
\begin{proof}
  % Définition_z de Phi
  Let $v \in V(\CalY)$. If $v \in V(Y)$ let $f \in \Atlas$ be the isometry
  defined on $B_Y(v,R)$. If $v$ is a print $P$ let $y$ and $f \in \Atlas$ be such that $P$ is
  a print associated to $y$ and $f$. Our goal is to show that the map
  \begin{equation*}
    \label{eq:locisomdeCalYversX}
    \phi_f \ : \
    \begin{cases}
      B_{\CalY}(v,\Rbdg) & \rightarrow \bdg,\\
      z \in Y & \mapsto qf(y),\\
      Q & \mapsto x \quad \text{where}\ \Print(x)=qf(Q),
    \end{cases}
  \end{equation*}
  is an isometry.
  
  % Bien déf + injective.
  By \cref{Prop:characterisationparprint}, it is a well defined map. Moreover,
  using the injectivity of $q$ and \cref{Prop:characterisationparprint} and
  \cref{eq:Partitonbdg} we get that $\phi_f$ is an injective map.
  
  % --------------------- Arêtes ---------------
  Now, recall that since $q$ is a quasi-isometry, two elements $q(x_1)$ and
  $q(x_2)$ joined by an edge in $\bdg$ might be at distance greater than $1$ in 
  $\myg$. If we want to prove that $\phi_f$ is surjective on
  $B_{\bdg}(\phi_f(v), \Rbdg)$ and preserves the distance, we have to show that
  there exists a radius $\rprint$ allowing us 
  to “reconstruct” all the edges of $B_{\bdg}(\phi_f(v),\Rbdg)$ in
  $B_{\CalY}(v,\Rbdg)$. Let $L, \varepsilon>0$ be
  such that $q$ is a $(L,\varepsilon)$-quasi-isometry. We distinguish three
  cases, represented in \cref{fig:preuvePhi}.\\
  If $\chi_1, \chi_2 \in im(q)$, then let $x_1, x_2 \in \myg$ such that
  $q(x_i)=\chi_i$. They verify $d_{\myg}(x_1,x_2) \leq
  Ld_{\bdg}(\chi_1,\chi_2) +\varepsilon$. This case is represented in
  \cref{fig:Phicas1}.\\ 
  If $\chi_1 \in im(q)$ and $\chi_2 \notin im(q)$, let $x_1= q^{-1}(\chi_1)$.
  For all $x_2 \in \myg$ such that $q(x_2) \in \Print(\chi_2)$, we have (see \cref{fig:Phicas2})
  \begin{equation*}
    d_{\bdg}(q(x_1), q(x_2)) \leq 1 + d_{\bdg}(\chi_1,\chi_2)  \ \Rightarrow
    d_{\myg}(x_1,x_2) \leq Ld_{\bdg}(\chi_1,\chi_2) + L+\varepsilon.
  \end{equation*}
  If $\chi_1, \chi_2 \notin im(q)$, let $x_i \in \myg$ such that $q(x_i) \in
  \Print(\chi_i)$ for $i=1,2$. Then (see \cref{fig:Phicas2})
  \begin{equation*}
    d_{\bdg}(q(x_1), q(x_2)) \leq 2 +d_{\bdg}(\chi_1,\chi_2)  \ \Rightarrow
    d_{\myg}(x_1,x_2) \leq Ld_{\bdg}(\chi_1,\chi_2) + 2L+\varepsilon.
  \end{equation*}
  % -----------------------------------------------------------------
  % La grande figure avec les trois cas
  % -----------------------------------------------------------------
  \begin{figure}[htbp]
    \begin{subfigure}[b]{0.32\textwidth}
      \centering
      \begin{tikzpicture}
        % Chi 2
        \node[draw, circle, fill=MFCB] (chi2) at (0,0) {};
        % q(x2)
        \node[draw, circle, fill=DeepSkyBlue4, thick] (blob) at (120:1) {};
        \node[anchor=south] (cchi2) at (120:1.1) {$\chi_2$};
        % Les autres voisins de chi2
        \node[draw, circle, fill=DeepSkyBlue4] (blob2) at (240:1) {};
        \node[draw, circle, fill=DeepSkyBlue4!20, dashed] (blob3) at (0:1) {};
        \draw (blob) -- (chi2) -- (blob2);
        \draw[dashed] (blob3) -- (chi2);
        \begin{scope}[shift={(240:1)}]
          \node[draw, circle, fill=MFCB] (geo) at (310:1) {};
          \node[draw, circle, fill=MFCB!20, dashed] (azertyu) at (190:1) {};
          \draw[dashed](azertyu) -- (blob2);
          \begin{scope}[shift={(310:1)}]
            \node[draw, circle, fill=DeepSkyBlue4, thick] (chi1) at (240:1) {};
            \node[anchor=north west] (lchi1) at (245:1.1) {$\chi_1$};
            \node[draw, circle, fill=DeepSkyBlue4!20, dashed] (az) at (0:1) {};
            \draw (chi1) -- (geo) -- (blob2);
            \node (blanc) at  (270:1.8) {};
          \end{scope}
          \draw[dashed](az) -- (geo);
        \end{scope}
      \end{tikzpicture}
      \caption{First case}
      \label{fig:Phicas1}
    \end{subfigure}
    \begin{subfigure}[b]{0.32\textwidth}
      \centering
      \begin{tikzpicture}
        % Chi 2
        \node[draw, circle, fill=MFCB, thick] (chi2) at (0,0) {};
        \node[anchor=north west] (lchi2) at (375:0.01) {$\chi_2$};
        % q(x2)
        \node[draw, circle, fill=DeepSkyBlue4] (blob) at (120:1) {};
        \node[anchor=east] (qx2) at (120:1.1) {$q(x_2)$};
        % Les autres voisins de chi2
        \node[draw, circle, fill=DeepSkyBlue4] (blob2) at (240:1) {};
        \node[draw, circle, fill=DeepSkyBlue4!20, dashed] (blob3) at (0:1) {};
        \draw (blob) -- (chi2) -- (blob2);
        \draw[dashed] (blob3) -- (chi2);
        \begin{scope}[shift={(240:1)}]
          \node[draw, circle, fill=MFCB] (geo) at (310:1) {};
          \node[draw, circle, fill=MFCB!20, dashed] (azertyu) at (190:1) {};
          \draw[dashed](azertyu) -- (blob2);
          \begin{scope}[shift={(310:1)}]
            \node[draw, circle, fill=DeepSkyBlue4, thick] (chi1) at (240:1) {};
            \node[anchor=north west] (lchi1) at (245:1.1) {$\chi_1$};
            \node[draw, circle, fill=DeepSkyBlue4!20, dashed] (az) at (0:1) {};
            \draw (chi1) -- (geo) -- (blob2);
            \node (blanc) at  (270:1.8) {};
          \end{scope}
          \draw[dashed](az) -- (geo);
          \node[anchor=south, color=black!50] (px2) at (70:2.5) {$\mathcal{P}(\chi_2)$};
          \draw[->, >=latex, black!30] (px2) to[bend right] (blob.north);
          \draw[->, >=latex, black!30] (px2) to[bend left] (blob3.north);
        \end{scope}
        \draw[->, >=latex, black!30] (px2) to[bend left] (355:1.4) to[bend left](blob2.east);
      \end{tikzpicture}
      \caption{Second case}
      \label{fig:Phicas2}
    \end{subfigure}
    % ----------------------------
    \begin{subfigure}[b]{0.32\textwidth}
      \centering
      \begin{tikzpicture}
        % Chi 2
        \node[draw, circle, fill=MFCB, thick] (chi2) at (0,0) {};
        \node[anchor=north west] (lchi2) at (375:0.01) {$\chi_2$};
        % q(x2)
        \node[draw, circle, fill=DeepSkyBlue4] (blob) at (120:1) {};
        \node[anchor=east] (qx2) at (120:1.1) {$q(x_2)$};
        % Les autres voisins de chi2
        \node[draw, circle, fill=DeepSkyBlue4] (blob2) at (240:1) {};
        \node[draw, circle, fill=DeepSkyBlue4!20, dashed] (blob3) at (0:1) {};
        \draw (blob) -- (chi2) -- (blob2);
        \draw[dashed] (blob3) -- (chi2);
        \begin{scope}[shift={(240:1)}]
          \node[draw, circle, fill=MFCB, thick] (geo) at (310:1) {};
          \node[anchor=north] (cchi1) at (310:1.15) {$\chi_1$};
          \node[draw, circle, fill=MFCB!20, dashed] (azertyu) at (190:1) {};
          \draw[dashed](azertyu) -- (blob2);
          \begin{scope}[shift={(310:1)}]
            \node[anchor=east, color=black!50] (px1) at (220:1.8) {$\mathcal{P}(\chi_1)$};
            \node[draw, circle, fill=DeepSkyBlue4] (chi1) at (240:1) {};
            \node[anchor=north west] (lchi1) at (245:1.1) {$q(x_1)$};
            \node[draw, circle, fill=DeepSkyBlue4!20, dashed] (az) at (0:1) {};
            \draw (chi1) -- (geo) -- (blob2);
            \draw[->, >=latex, black!30] (px1) to[bend right] (270:1.8) to[bend right] (az.south);
          \end{scope}
          \draw[->, >=latex, black!30] (px1) to[bend left] (chi1.west);
          \draw[->, >=latex, black!30] (px1) to[bend left] (blob2.west);
          \draw[dashed](az) -- (geo);
          \node[anchor=south, color=black!50] (px2) at (70:2.5) {$\mathcal{P}(\chi_2)$};
          \draw[->, >=latex, black!30] (px2) to[bend right] (blob.north);
          \draw[->, >=latex, black!30] (px2) to[bend left] (blob3.north);
        \end{scope}
        \draw[->, >=latex, black!30] (px2) to[bend left] (355:1.4) to[bend left](blob2.east);
      \end{tikzpicture}
      \caption{Third case}
      \label{fig:Phicas3}
    \end{subfigure}
    \caption{The three cases ($im(q)$ is represented by the blue vertices)}
    \label{fig:preuvePhi}
  \end{figure}
  % -----------------------------------------------------------------
  %  figure avec les trois cas - FIN
  % -----------------------------------------------------------------
  Hence, assume $\rprint > L\Rbdg + 2L+\varepsilon$ and let us show that $\phi_f$
  is an isometry.

  Let $\chi \in B_{\bdg}(\phi_f(v),\Rbdg)$, by choice of $\rprint$ either $\chi
  \in im(q)$ and then there exists $z \in B_Y(y,\rprint)$ such that $qf(z)=\chi$
  or $\chi \notin im(q)$ and then there exists $P \subset B_{\mY}(y,\rprint)$
  such that $qf(P)=\Print(\chi)$.
  Hence, in both cases $\chi \in im(\phi_f)$ and thus, $\phi_f$ is a bijection
  from $B_{\CalY}(v,\Rbdg)$ to $B_{\bdg}(\phi_f(v),\Rbdg)$.
  Now take $v_1, v_2 $ in $B_{\CalY}(v,\Rbdg)$ at distance $l$ in $\CalY$ and
  let $(w_0=v_1, w_1, \ldots, w_l=v_2)$ be a geodesic in $\CalY$. By definition
  of $\CalY$ and choice of $\rprint$, for all $i \in \{0, \ldots,l-1\}$ if there
  is an edge between 
  $w_i$ and $w_{i+1}$, then $d(\phi_f(w_i),\phi_f( w_{i+1}))=1$. Hence
  $d_{\bdg}(\phi_f(v_1), \phi_f(v_2)) \leq l$. 
  To get the reversed inequality, take $\chi_1, \chi_2$ in 
  $B_{\bdg}(\phi_f(v),\Rbdg)$. Since $\phi_f$ is bijective there exists $v_0,
  \ldots, v_l$ in $\CalY$ such that $(\phi_f(v_0), \ldots, \phi_f(v_l))$ is a
  geodesic between $\chi_1$ and $\chi_2$. Again, by definition of $\CalY$ and
  choice of $\rprint$, an edge between $\phi_f(v_i)$ and $\phi_f(v_{i+1})$ gives
  an edge between $v_i$ and $v_{i+1}$ in $\CalY$ and thus $d_{\CalY}(v_1, v_2) \leq
  l$.
  
  Hence, if $\rprint >L\Rbdg + 2L+\varepsilon$ then $\phi_f$ is an isometry.
\end{proof}
The LG-rigidity of the building will give us a covering from $\bdg$ to $\calY$.
In order to obtain an \emph{isometry} we need to prove (by \cref{Prop:15DlST2})
that $\calY$ is simply connected at the same scale as $\bdg$. 

\begin{Lmm}\label{Lmm:CalY3sco}
  If $\Rbdg$ (and hence $R$) is large enough, then $\CalY$ is simply connected
  at scale $3$. 
\end{Lmm}
We first prove that $\CalY$ is quasi-isometric to $\mY$ and use it to show that
$\CalY$ is simply connected at some scale $k^\prime$. We conclude using the
contractibility of the building and the fact
that $\CalY$ is locally the same as the building. But before looking at the
detail of the proof, let us make a remark. 
\begin{Rq}\label{Rq:DistancePrintPi}
  Let $P$ be a print associated to some $z\in \mY$ and $f\in \Atlas$ and
  let $y\in P$. If $x$ is the source of $qf(P)$, then 
  $d_{\CalY}(P,y)=d_{\bdg}(x,qf(y))=1$. 
\end{Rq}
\begin{proof}[Proof of \cref{Lmm:CalY3sco}]
  Let us show that $\CalY$ is quasi-isometric to $\mY$. Define $\pi:$ $\CalY
  \rightarrow \mY$ such that if $y\in V(Y)$ then $\pi(y)=y$ and if $P$ is a
  print then $\pi(P)=y$ for some $y\in P$ arbitrarily chosen. Let $(v_0, \ldots,
  v_m)$ be a geodesic in $\CalY$ and for all $i\in \{0,\ldots,m\}$ define
  $y_i:=\pi(v_i)$ and $f_i$ to be the isometry of $\Atlas$ associated to
  $y_i$. Using that $q$ is a $(L,\varepsilon)$-quasi-isometry, we get
  \begin{align*}
    d_{\mY}\left(\pi(v_0), \pi(v_m)\right)=d_{\mY}(y_0,y_m)
    &\leq \sum^m_{i=0}d_{\mY}\left( y_i,y_{i+1}\right),\\
    &\leq \sum^m_{i=0}\left[Ld_{\bdg}\left(qf_i(y_i),qf_i(y_{i+1})\right)
      +\varepsilon  \right].
  \end{align*}
  Now let $i\in \{0,\ldots,m\}$. If $v_i$ is a print, denote by $x_i$ the source
  of $qf(v_i)$ and if $v_i$ belongs to the copy of $V(\mY)$ contained in
  $\CalY$ let $x_i:=qf_i\pi(v_i)$. Then $d_{\CalY}\left(v_i,v_{i+1}\right)=
  d_{\bdg}(x_i,x_{i+1})$ for all $i$. Thus, using \cref{Rq:DistancePrintPi}, we get
  \begin{align*}
    d_{\bdg} \left( qf_i(y_i),qf_i(y_{i+1})\right)
    &\leq d_{\bdg}\left( qf_i(y_i),x_i\right)+ d_{\bdg}(x_i,x_{i+1})
      + d_{\bdg}\left(qf_i(y_{i+1}),x_{i+1}\right),\\
    &\leq 2+ d_{\bdg}(x_i,x_{i+1})=2+d_{\CalY}\left(v_i,v_{i+1}\right).
  \end{align*}
  Since $d_{\CalY}\left(v_i,v_{i+1}\right)=1$, we obtain
  \begin{align*}
    d_{\mY}\left(\pi(v_0), \pi(v_m)\right)=d_{\mY}(y_0,y_m)
    &\leq \sum^m_{i=0}\left[L2+Ld_{\CalY}\left(v_i,v_{i+1}\right)+\varepsilon  \right],\\
    &=(3L+\varepsilon)m=(3L+\varepsilon)d_{\CalY}\left(v_0,v_{m}\right).
  \end{align*} 
  Now let $v,v^\prime \in \CalY$ and let $(\pi(v)=z_0,\ldots,\pi(v^\prime)=z_l)$
  be a geodesic in $\mY$. For all $i\in \{0,\ldots,l\}$ take $f^{\prime}_i\in \Atlas$ the
  isometry associated to $z_i$. Then
  \begin{equation*}
    d_{\CalY}(v,v^\prime)
    \leq d_{\CalY}(v,z_0)+\sum^{l-1}_{i=0}d_{\CalY}(z_i,z_{i+1})
   +d_{\CalY}(z_l,v^\prime).
  \end{equation*}
  But by \cref{Rq:DistancePrintPi} if $v$ (resp. $v^\prime$) is a print then
  $d_{\CalY}(v,z_0)=1$ (resp. $d_{\CalY}(v^\prime,z_l)=1$). And if $v$ (resp.
  $v^\prime$) belongs to $V(Y)$ then $v=z_0$ (resp. $v^\prime=z_l$). Thus both
  $d_{\CalY}(v,z_0)$ and $d_{\CalY}(v^\prime,z_l)$ are always smaller than $1$.
  Hence,
  \begin{align*}
    d_{\CalY}(v,v^\prime)
    \leq 2 + \sum^{l-1}_{i=0}d_{\CalY}(z_i,z_{i+1})
    &= 2 + \sum^{l-1}_{i=0}d_{\bdg}\left(qf^{\prime}_i(z_i),qf^{\prime}(z_{i+1})\right),\\
    &\leq 2+ \sum^{l-1}_{i=0}\left[Ld_{\mY}(z_i,z_{i+1})+\varepsilon\right],\\
    &= 2+ (L+\varepsilon)l=2+(L+\varepsilon)d_{\mY}(\pi(v),\pi(v^\prime)).
  \end{align*}
  Thus $\pi$ is a quasi-isometry between $\CalY$ and $\mY$. Hence \cref{Prop:Th22DlST2} implies
  that there exists $k^\prime \in \bN^{*}$ such that $\CalY$ is simply-connected
  at scale $k^\prime$. 
  
  Finally, let $\ell$ be loop in $\CalY$ of length less than 
  $k^\prime$. If $\Rbdg$ is large enough then $\ell$ is contained in some
  ball $B$ in $\calY$. By \cref{Lmm:calYRlocX} there exists a local isometry $\phi$ from
  $B$ to some ball $\mathcal{B}$ in $\bdg$. But $\phi(\ell)$ is
  contractible inside its convex hull, by \cref{Claim:ContractibilityofBdg}. In
  particular it is simply-connected. Since $\bdg$ is $3$-simply-connected and if
  $\Rbdg$ is large enough, the
  convex hull of $\phi(\ell)$ is contained in the complex obtained by gluing
  triangles on all the loops of length $3$ in $\mathcal{B}$. Which, by local
  isometry with $B$, proves the wanted assertion.
\end{proof}
Thanks to the previous lemma, we can now use the rigidity of the Bruhat-Tits building.

\begin{Prop}\label{Prop:CalYisombdg}
  If ${\Rbdg}$ (and hence $R$) is large enough, then $\CalY$ is isometric to $\bdg$.
\end{Prop}

\begin{proof} Recall that we have $R > \ratlas >3\rprint> \rprint > 3\Rbdg + 2L
  + \varepsilon> \Rbdg$. \\
  By \cref{Th:LGRImmeuble}, the building $\bdg$ is LG-rigid.
  Moreover, since its isometry group is transitive 
  \cref{Prop:15DlST2} gives us the existence of some radius 
  $R_{sc}>0$ such that every graph which is $3$-simply connected and
  $R_{sc}$-locally $\bdg$ is isometric to $\bdg$.\\
  By definition of the edges on $\CalY$, this graph is simply connected at scale
  $3$. Taking $\rprint$ (and hence $R$) large enough so that ${\Rbdg} \geq R_{sc}$
  the preceding paragraph combined with \cref{Lmm:calYRlocX} give us the
  existence of an isometry between $\bdg$ and $\CalY$.
\end{proof}
\subsection{Change of local map, change of global isometry}
Let $y \in \mY$ and $f_y \in \Atlas$ be the isometry defined on $B(y,R)$. Let  
\begin{equation}
  \label{eq:locisomdeCalYversX0}
  \phi_y \ : \
  \begin{cases}
    B_{\CalY}(y,\Rbdg) & \rightarrow \bdg\\
    z \in Y & \mapsto qf_y(z)\\
    Q & \mapsto x \qquad \quad\text{where} \ \Print(x)=qf_y(Q).
  \end{cases}
\end{equation}
\begin{Lmm} \label{Lmm:transitionmapdecalY}
  Let $y$ and $z$ be neighbours in $Y$ and $a \in H_0$ such that $f_yf^{-1}_z$
  coincide with $a$ on $B_{\myg}(f(z), \ratlas)$. If $\Rbdg$ is large enough, then
  $\phi_y\phi^{-1}_z$ coincide with $\rho(a)$ on $B_{\bdg}(\phi_z(z),2)$.
\end{Lmm}
\begin{proof}
  Let $y$ and $z$ be neighbours in $Y$ and $a \in H_0$ such that $f_yf^{-1}_z$
  coincide with $a$ on $B_{\myg}(f(z), \ratlas)$. If $\Rbdg$ (and hence $R$) is large
  enough, then $B_{\CalY}(z,2)$ is contained in $B_{\CalY}(y,\Rbdg)$. Thus,
  $\phi_y\phi^{-1}_z$ is well defined on $B_{\bdg}(\phi_z(z),2)$.

  Let $v \in B_{\CalY}(z,2)$. If $v \in V(Y)$, then
  \begin{equation*}
    \phi_y(v) = qf_y(v) = qaf_z(v)=\rho(a)qf_z(v)=\rho(a)\phi_z(v). 
  \end{equation*}
  If $v=P$ with $P\subset Y$ a print, then
  \begin{equation*}
    \Print \big(\phi_y(v) \big)= qf_y(P) = qaf_z(P)=\rho(a)qf_z(P)=
    \Print \big(\rho(a)\phi_z(v) \big),
  \end{equation*}
  Thus $\phi_y(v)=\rho(a)\phi_z(v)$, since the print determines the vertex.
  Hence the result.
\end{proof}

Now let $\rX >0$. If $\Rbdg$ is large enough then, by
SLG-rigidity\label{firstmentionofSLGX} of $\bdg$ there exists an isometry
$\iota_y$ from $\CalY$ to $\bdg$ that coincides with $\phi_y$ on $B(y,\rX)$.
Thus, the lemma above allows us to work with a set of isometries from $\CalY$ to $\bdg$
that differs only by a multiplication by an element of $PSL_n(\bK)$.
\begin{Lmm}\label{Claim:iotapreservepartition}
  If $y$ and $z$ belong to $Y$ and $\Rbdg$ is large enough, then
  $\iota_y\iota^{-1}_z \in PSL_n(\bK)$. Hence for all $y \in \mY$, the
  isometry $\iota_y$ sends the copy of  
  $V(\mY)$ contained in $\CalY$ to $im(q)$ and sends prints contained in $\CalY$ to
  vertices in $\bdg \backslash im(q)$.
\end{Lmm}
\begin{proof}
  Let $y$ and $z$ be neighbours in $Y$. Since $\iota_y \iota^{-1}_z$ is an
  isometry of $\bdg$ it permutes the $PSL_n(\bK)$-orbits. Recall that
  $\iota_y$ coincides with $\phi_y$ on $B(y,\rX)$. 
  Hence, if $\rX$ (and
  hence $R$) is large enough, then $B_{\CalY}(z,2)$ is contained in
  $B_{\CalY}(y,\rX)$, thus 
  \begin{equation*}
    {\left(\iota_y\iota^{-1}_z  \right)}_{|B_{\bdg}(\iota_z(z),2)} = \phi_y\phi^{-1}_z. 
  \end{equation*}
  But $\phi_y\phi^{-1}_z$ coincides with an element of $PSL_n(\bK)$ on
  $B_{\bdg}(\phi_z(z),2)$, by \cref{Lmm:transitionmapdecalY}.
  Hence $\iota_y\iota^{-1}_z$ restricted to a ball of radius $2$ preserves the
  $PSL_n(\bK)$-orbits. Since such a ball contains a vertex of each type, it
  implies that $\iota_y\iota^{-1}_z$ preserves the $PSL_n(\bK)$-orbits and thus
  belongs to $PSL_n(\bK)$.

  Now take $y$ and $z$ in $Y$ (not necessarily neighbours). %Denote by
  Let $(y_0=y,y_1,\ldots, y_l=z)$ be a geodesic in $\mY$. By the preceding paragraph, 
  there exists a sequence $\alpha_1$, …, $\alpha_{l}$ of elements in
  $PSL_n(\bK)$ such that
  \begin{equation*}
    \forall i \in \{1, \ldots, l\} \quad \iota_{y_i}\iota^{-1}_{y_{i-1}} = \alpha_i.
  \end{equation*}
  Thus, recalling that $z=y_l$ and $y=y_0$, we get $\iota_{z}=\alpha_l \cdots
  \alpha _1 \iota_{y}$. Which proves the first assertion of the lemma.

  Let us now prove the second part of the lemma. Let $y \in Y$ and $v \in
  \CalY$. There exists $z \in \mY$ such that $v \in B_{\CalY}(z,2)$, and
  using the paragraph above, there exists $\alpha \in PSL_n(\bK)$ such that
  $\iota_{y}=\alpha \iota_{z}$. In particular, since $v$ belongs to $B_{\CalY}(z,\Rbdg)$,
  \begin{equation*}
    \iota_y(v)= \alpha \iota_z(v)
    =  \alpha\phi_z(v).
  \end{equation*}
  By definition of $\phi_z$, if $v \in V(Y)$ then $\phi_z(v)$ belongs to $im(q)$
  and if $v =P$ with $P\subset Y$ a print, then $\phi_z(v)$ belongs to $\bdg
  \backslash im(q)$. This finishes the proof of the lemma.
\end{proof}

Now we have all the tools we need to prove the isometry between $\mY$ and $\myg$.
\subsection{\texorpdfstring{Isometry from $\mY$ to $\myg$}{Isometry between Y
    and X}}
\label{sec:isommYmyg}
Let $\kappa$ be the natural injection of $Y$ in $ \calY_Z$ and $\iota$ an
isometry given by \cref{Prop:CalYisombdg}. With the objects constructed so far
we get the diagram in \cref{fig:schemadebutpartiequatre}.

\begin{figure}[htbp]
  \centering
  \begin{subfigure}[b]{0.42\textwidth}
    \centering
    \begin{tikzpicture}
      % \myg and Y rloc myg
      \node[Leturquoise] (myg) at (0,2) {$\myg$};
      \node[Leturquoise] (mY) at (0,0) {$\mY$};
      \node[Leturquoise] (y) at (0,-0.5) {$y$};
      \draw[<->, >=latex, dashed] (myg.south) -- (mY.north) node[left, midway] {$R$-loc};
      %Immmeuble and QI 
      \begin{scope}[shift={(1,0)}]
        \node[MFCB] (CalY) at (2,0) {$\calY$};
        \node[MFCB] (ky) at (2,-0.5) {$\kappa(y)$};
        \node[MFCB] (bdg) at (2,2) {$\bdg$};
      \end{scope}
      \draw[{Hooks[right]}->, >=latex] (myg.east) -- (bdg.west)
      node[above, midway] {$q$};
      \draw[<-, >=latex] (bdg.south) -- (CalY.north)
      node[right, midway] {$\iota$, isom.};
      \draw[{Hooks[right]}->, >=latex] (mY.east) -- (CalY.west)
      node[below, midway] {$\kappa$};
      \draw[{Bar}->, >=latex] (y.east) -- (ky.west);
      \node (edvbjkldkfr) at (0,-1) {};
    \end{tikzpicture}  
    \caption{Maps between graphs}
    \label{fig:mapsandgraphs}
  \end{subfigure}
  \hfill 
  \begin{subfigure}[b]{0.42\textwidth}
    \centering
    \begin{tikzpicture}
      % Bloc isometries
      \begin{scope}[shift={(0,1.5)}]
        \begin{scope}[shift={(1,0)}]
          \node[MFCB] (IsX) at (2,3.5) {$\Is(X)$};
          \node[MFCB] (rhoa) at (2,3) {$\rho(a)$};
          \node[MFCB] (PSL) at (2,2.5) {$PSL_n(\bK)$};
        \end{scope}
        \node[Leturquoise] (IsGamma) at (0,3.5) {$\Is(\myg)$};
        \node[Leturquoise] (a) at (0,3) {$a$};
        \draw[{Hooks[right]}->, >=latex] (IsGamma) -- (IsX) node[above, midway] {$\rho$}; 
        \draw[{Bar}->] (a) -- (rhoa); 
        \node[Leturquoise] (H0) at (0,2.5) {$H_0$};
        \draw[->, >=latex] (H0) -- (PSL);
        \node (edvbjkldkfr) at (0,1) {};
      \end{scope}
    \end{tikzpicture}
    \caption{Relations between groups}
    \label{fig:Relgroups}
  \end{subfigure}
  \caption{Relations between the different graphs and groups}
  \label{fig:schemadebutpartiequatre}
\end{figure}

The aim of this section is to prove the following result.

\begin{Prop}\label{Prop:Yisommyg}
  For $\Rbdg$ large enough, the graphs $\mY$ and $\myg$ are isometric.
\end{Prop}
Let us discuss the strategy of the proof. Using the preceding section, we chose
an isometry $\iota$ from $\CalY$ to $\bdg$ that coincides with a $\phi_y$ on a
small ball. Then,
we show that $\kappa \iota q^{-1}$
is locally an isometry, \textit{viz.} there exists a radius $\rY$ 
such that ${q^{-1}} \iota \kappa$ restricted to any ball of radius $\rY$
preserves the distance. We conclude by showing that it forces $\kappa \iota
q^{-1}$ to be an isometry.

\begin{proof}[Proof of \cref{Prop:Yisommyg}] By
  \cref{Claim:iotapreservepartition}, for any $y\in Y$ the map $q^{-1} \iota_y
  \kappa $ is well defined.
  Now fix $y_0 \in Y$ and consider $\iota:=\iota_{y_0}$. We want to prove that
  $q^{-1} \iota \kappa$ restricted to small balls preserves the distance. Then
  we will show that it is an isometry from $\mY$ to $\myg$.
  \begin{Claim} \label{Claim:localeisometrie}
    Let $y \in Y$ and $\rY \geq 1$. If $R$ is large enough, then $q^{-1}\iota \kappa$
    restricted to $ B_{\mY}(y,\rY)$ preserves the distance. 
  \end{Claim}
  \begin{proof}[Proof of the claim]
    Let $\rY \geq 1$ and recall that we have $ R > \ratlas >3\rprint> \rprint >
    3\Rbdg + 2L + \varepsilon> \Rbdg > \rX$.
    Let $y \in \mY$ and recall that $L$ and $\varepsilon$ are
    constants such that $q$ is a $(L,\varepsilon)$-quasi-isometry. If $\rX
    \geq L\rY + \varepsilon$ (and hence if $R$ is large enough) then
    $\kappa(B_{\mY}(y,\rY))$ is included in 
    $B_{\CalY}(y,\rX)$.  Indeed if $z\in B_{\mY}(y,\rY)$ then
    \begin{equation*}
      d_{\bdg} \Big(qf_y(y),qf_y(z) \Big) \leq  L d_{_\myg} \Big(f_y(y),f_y(z) \Big)
      + \varepsilon = Ld_{\mY}(y,z) + \varepsilon \leq L \rY + \varepsilon \leq \rX. % \Rbdg.
    \end{equation*}
    Thus $\phi_y(\kappa(z))=qf_y(z)$ and
    \begin{align*}
      d_{\CalY}\big( \kappa(y),\kappa(z) \big)
      = d_{\bdg}\left(\phi_y\big(\kappa(y)\big),\phi_y\big(\kappa(z)\big) \right)
      &=  d_{\bdg}\Big(qf_y(y),qf_y(z) \Big)
      \leq \Rbdg.
    \end{align*}

    Now, recall that
    $H_0=\rho^{-1}PSL_n(\bK)$. Then, by \cref{Claim:iotapreservepartition}
    there exists $a_y \in H_0$ such that $\iota_y \iota^{-1}=\rho(a_y)$. Hence,
    using the equivariance of $q$ we get that for all $z_1$ and $z_2$ in $B_{\mY}(y,\rY)$ 
    \begin{align*}
      d_{\myg} \Big(q^{-1}\iota \kappa(z_1),q^{-1}\iota \kappa(z_2)\Big)
      &=d_{\myg} \Big(a_yq^{-1}\iota \kappa(z_1),a_yq^{-1}\iota \kappa(z_2)\Big)\\
      & =d_{\myg} \Big(q^{-1}\rho(a_y)\iota \kappa(z_1),q^{-1}\rho(a_y)\iota \kappa(z_2)\Big)\\
      & =d_{\myg} \Big(q^{-1}\iota_y \kappa(z_1),q^{-1}\iota_y \kappa(z_2)\Big).
    \end{align*}
    But $z_1$ and $ z_2 $ belong to $ B_{\mY}(y,\rY)$, hence for $i=1,2$ we have $\iota_y\kappa(z_i)
    = qf_y(z_i)$. Thus,
    \begin{align*}
      d_{\myg} \Big(q^{-1}\iota \kappa(z_1),q^{-1}\iota \kappa(z_2)\Big)
      &=d_{\myg} \Big(q^{-1} qf_y(z_1),q^{-1} qf_y(z_2)\Big)\\
      &=d_{\myg} \Big(f_y(z_1),f_y(z_2)\Big)
        = d_{\mY}(z_1,z_2).
    \end{align*}
    Thus $q^{-1}\iota \kappa$ restricted to $ B_{\mY}(y,\rY)$ preserves the distance.
  \end{proof}
  Let’s show that the claim forces $q^{-1}\iota \kappa$ to
  be an isometry from $\mY$ to $\myg$. Take $\rY\geq 2$ and let $y, y^\prime \in
  \mY$ and $(y_0=y, y_1, \ldots, y_l=y^\prime)$ be a geodesic in $\mY$. Since
  for all $i$ the vertices $y_i$ and $y_{i+1}$ are adjacent, then
  \cref{Claim:localeisometrie} implies that $d_{\myg}(q^{-1}\iota \kappa
  (y_i),q^{-1}\iota \kappa (y_{i+1}))=1$. Hence
  \begin{align*}
    d_{\myg}\Big(q^{-1}\iota \kappa (y),q^{-1}\iota \kappa (y^\prime)\Big)
    \leq \sum^{l-1}_{i=0} d_{\myg}\Big(q^{-1}\iota \kappa (y_i),q^{-1}\iota \kappa (y_{i+1})\Big)
    =l.
  \end{align*}
  Moreover, if $(x_0=q^{-1}\iota \kappa(y), x_1, \ldots, x_m=q^{-1}\iota
  \kappa(y^\prime))$ is a geodesic in $\myg$, then by bijectivity of $q^{-1}\iota
  \kappa$ there exists $z_i \in \mY$ such that $q^{-1}\iota
  \kappa(z_i) = x_i$ for all $i$ in $\{1, \ldots, m-1\}$. Denote $z_0=y$ and $z_m=y^\prime$. Since
  for all $i$ the vertices $x_i$ and $x_{i+1}$ are adjacent, then
  \cref{Claim:localeisometrie} implies that
  $d_{\myg}(z_i,z_{i+1})=d_{\bdg}(q^{-1}\iota \kappa (z_i), q^{-1}\iota \kappa (z_{i+1}))$. Thus
  \begin{align*}
    d_{\mY}(y,y^\prime) \leq \sum^{m-1}_{i=0} d_{\mY}(z_i,z_{i+1})
    &= \sum^{m-1}_{i=0} d_{\myg}\Big(q^{-1}\iota \kappa (z_i),q^{-1}\iota \kappa (z_{i+1})\Big),\\
    &= \sum^{m-1}_{i=0} d_{\myg}(x_i, x_{i+1})
    =m.
  \end{align*}
\end{proof}
We conclude by the proof of \cref{Th:BigRigidity}.
\begin{proof}[Proof of {\cref{Th:BigRigidity}}] Let $n\neq 3$ and $\myg$
  verifying the hypothesis of \cref{Th:BigRigidity}.
  
  If $n=2$ then $\bdg$ is the $(p+1)$-regular tree, thus by \cref{Ex:QuasiArbre}
  if $\myg$ is quasi-isometric to $\bdg$ then $\myg$ is LG-rigid.

  If $n\geq 4$, let $k\in \bN$ such that $\myg$ is simply connected at scale
  $k$. Then by \cref{Prop:Yisommyg} for $R$ large enough, any $k$-simply-connected graph
  $\mY$ being $R$-locally the same as $\myg$ is isometric to $\myg$. Thus $\myg$
  is LG-rigid.
  Finally for any $n\neq 3$, since $\myg$ is assumed transitive it is actually
  SLG-rigid by \cref{Prop:LGRSLGR}.  
\end{proof}
% Section 4
\section{\texorpdfstring{Application to $p$-adic lattices}{Application to
    p-adic lattices}}
In this section we prove \cref{Th:LGRReseaux} which we recall below. 
\begin{Cor} \label{Cor:LGRReseaux} Let $n\neq 3$ and $\bK$ be a non-Archimedean
  skew field of characteristic zero.\\
 The torsion-free lattices of $SL_n(\bK)$ are SLG-rigid. 
\end{Cor}

Let $n \neq 3$, let $\bK$ be a non-Archimedean
skew field of characteristic zero and $\Gamma\leq SL_n(\bK)$ be a lattice
without torsion. Denote by $(\Gamma, S)$ one of its Cayley graphs. Recall that
any lattice in $SL_n(\bK)$ is uniform (\textit{i.e.} cocompact).
\subsection{Quasi-isometry between the lattice and the building}
To show the corollary, we first check that the lattice is quasi-isometric to the
building. Then, using a famous result of Kleiner and Leeb we show that the
isometry group of the lattice acts on the building and that the quasi-isometry
can be chosen to be equivariant under this action. 
% ===================================
\begin{Lmm}\label{Lmm:GammaQIX}
  Let $\Lambda$ be a lattice of $SL_n(\bK)$.
  Then $\Lambda$ is quasi-isometric to $\bdg$.
\end{Lmm}

\begin{proof} First, recall that any lattice in $SL_n(\bK)$ is uniform,
  \textit{viz.} cocompact (see for example \cite{BenoistQuint}).
  
  Since $\Lambda$ is a lattice of $SL_n(\bK)$, there is a natural action on the 
  Bruhat-Tits building induced by the action of $PSL_n(\bK)$. Moreover, since
  $\Lambda$ is cocompact and the $PSL_n(\bK)$ action has exactly $n$ orbits,
  the $\Lambda$ action is also cocompact. Hence by the Svarc-Milnor's lemma
  $\Lambda$ is quasi-isometric to $X$.
\end{proof}

By a result of Kleiner and Leeb \cite{KL} and Cornulier \cite[Theorem
3.B.1]{CorKL} applied to our lattice $\Gamma$, 
this quasi-isometry implies the existence of a homomorphism from $\Is(\Gamma,S)$
to $\Is(X)$ and a quasi-isometry from $(\Gamma,S)$ to $X$ which is
$\Is(\Gamma,S)$-equivariant. 
Since $\Gamma$ is assumed to be torsion-free, we can refine the informations about
these two applications.

\begin{Lmm}
  \label{Lmm:KLaffine}
  Let $\Lambda$ be a lattice of $SL_n(\bK)$ and $T$ a symmetric generating
  set. If $\Lambda$ is torsion-free, then 
  there exists an injective homomorphism
  \begin{equation*}
    \rho \ : \ \Is(\Lambda,T) \ \rightarrow \Is(X),  
  \end{equation*}
  and an injective quasi-isometry which is $\Is(\Lambda,T)$-equivariant
  \begin{equation*}
    q \ : \ (\Lambda,T) \rightarrow X.
  \end{equation*}
\end{Lmm}

\begin{proof}
  Since we assumed that $\Lambda$ has no torsion element, by \cref{DlSTProp62}
  the isometry group of $(\Lambda,T)$ contains no non-trivial compact
  normal subgroup. Hence the morphism $\rho$ given by Kleiner-Leeb’s theorem is
  injective.

  Assume that there exist $\lambda_1, \lambda_2 \in \Lambda$
  such that $\lambda_1 \neq \lambda_2$ and $q(\lambda_1)=q(\lambda_2)$. 
  Then, the equivariance of $q$ implies that
  \begin{equation*}
    q\left(\left\{ {\left(\lambda_1\lambda^{-1}_2  \right)}^n\
        : \ n \in \bN\right\}\right)
    =\left\{q \left(e \right)\right\},
  \end{equation*}
  which contradicts the fact that $q$ is a quasi-isometry.
\end{proof}

\subsection{Relation between the isometry groups}
To apply \cref{Th:BigRigidity}, we still need to check that $\Is(\Gamma,S)$ is of
finite index in $\Is(\bdg)$. As stated in the lemma below, this is not always the
case: the lattice’s isometry group can also be discrete. But as we will see in
\cref{Preuverigiditedureseau}
we will be able to prove the rigidity of the lattice in that case too.

\begin{Lmm} \label{Lmm:Decomposition}
  Using the previous notations,
  \begin{itemize}
  \item Either $\Is(\Gamma,S)$ is discrete.
  \item Or $\Is(\Gamma,S)$ is of finite index in $\Is(X)$ and contains $PSL_n(\bK)$.
  \end{itemize}
\end{Lmm}
Before proving this lemma, let us recall a useful consequence of a theorem of Benoist and Quint.
The original and more general statement can be found in \cite[Corollary
4.5]{BenoistQuint}. 
\begin{Prop}[Benoist, Quint \cite{BenoistQuint}]
 Let $G$ be $p$-adic Lie group and $H$ be a finite covolume closed subgroup
 of $G$, with Lie algebra $\mathfrak{h}$. If $G$ has no proper cocompact normal
 subgroup, then $G$ normalizes $\mathfrak{h}$.
\end{Prop}

\begin{proof}[Proof of \cref{Lmm:Decomposition}]  
  Let $G=PSL_n(\bK)$ and  $H=\Is(\Gamma,S)\cap G$
  and note $\mathfrak{h}=:Lie(H)$ and $\mathfrak{g}:=Lie(G)$ their respective Lie
  algebras. Since $\Gamma$ is a lattice in $SL_n(\bK)$, we get that
  $\rho(\Gamma) \cap PSL_n(\bK)$ is a lattice in $PSL_n(\bK)$. Hence $H$
  contains the uniform lattice $\rho(\Gamma)\cap G$ of $G$, thus $H$ has finite
  covolume in $PSL_n(\bK)$.\\
  If $\bK$ is a non-Archimedean
  local skew field of characteristic zero then it is an extension of $\bQ_p$ for
  some prime $p$ (see for example \cite[Section 1]{DlST2016}). In particular $G$
  is a $p$-adic Lie group. Thus the above
  property applied to $G$ and $H$ implies that $G$ normalises $\mathfrak{h}$,
  in other words $\mathfrak{h}$ is an ideal of $\mathfrak{g}$. 
  Since $\mathfrak{g}$ is simple, we get that $\mathfrak{h}$ is either trivial
  or the full Lie algebra $\mathfrak{g}$.
  If $\Is(\Gamma,S)$ isn’t discrete, then it is a closed subgroup of
  $\Is(X)$. Hence $H$ is a closed subgroup of $G$ and its Lie algebra is
  non-trivial. By the previous point it can
  only be $\mathfrak{g}$. Hence, it implies that $H$
  is an open subgroup of $G$. Since it is also 
  cocompact, it is necessarily of finite index in $G$. Thus, we get that
  $\rho \left(\Is(\Gamma,S)\right)$ is of finite index in $\Is(X)$.
  \smallskip

  Let's show that $PSL_n(\bK) \leq \rho \left(\Is\left( \Gamma,S\right)\right)$.
  First assume that $\rho(\Is(\Gamma,S))$ is strictly contained in $PSL_n(\bK)$.
  Since these two groups are of finite index in $\Is(X)$, we get that
  $\rho\left(\Is(\Gamma,S) \right)$ is of finite index in $PSL_n(\bK)$.
  But then the core:
  \begin{equation*}
    \bigcap_{g\in PSL_n} g \cdot \rho\left(\Is(\Gamma,S) \right)\cdot g^{-1}
  \end{equation*}
  of $\rho\left(\Is(\Gamma,S)\right) $ is itself of finite index in $PSL_n(\bK)$ (and
  different from $PSL_n(\bK)$), which contradicts the simplicity of
  $PSL_n(\bK)$.
  
  Now, let's go back to the general case. Assume that $PSL_n(\bK)$ isn’t
  included in $\rho(\Is(\Gamma,S))$ and remark that:
  \begin{equation*}
    \mathfrak{h} = Lie \left( \Is(X) \right) = Lie \left( PSL_n(\bK)\right).
  \end{equation*}
  In particular $\rho(\Is(\Gamma,S))$ is “locally” $PSL_n(\bK)$ so,
  up to apply what precedes to an open set centered on $e_{\Gamma}$ sufficiently
  small of $\rho(\Is(\Gamma,S))$, we obtain a contradiction.\\
  Hence $PSL_n(\bK)$ is contained in $ \rho(\Is(\Gamma,S))$.
\end{proof}
\subsection{\texorpdfstring{Rigidity of $p$-adic lattices}{Rigidity of p-adic
    lattices}}
We conclude by the proof of \cref{Cor:LGRReseaux}.
\label{Preuverigiditedureseau}
\begin{proof}[Proof of \cref{Cor:LGRReseaux}] Let $n\neq 3$ and $p$ be a prime. %\\
  Let $\Gamma$ be a torsion-free lattice of $PSL_n(\bK)$ and $S$ be a
  symmetric generating part.\\
  If $n=2$, then $\bdg$ is the $(p+1)$-regular tree. Since by
  \cref{Lmm:GammaQIX}, the graph $(\Gamma,S)$ is
  quasi-isometric to $\bdg$, \cref{Ex:QuasiArbre} implies that $(\Gamma,S)$ is
  LG-rigid.

  Assume now that $n>3$. If $\Is(\Gamma,S)$ is discrete the LG-rigidity of the lattice is given by 
  \cref{Th:DlSTdiscrete}. \\
  If $\Is(\Gamma,S)$ is non-discrete, then 
  by \cref{Lmm:Decomposition} 
  it has finite index in $\Is(X)$ and in this case the hypothesis
  of \cref{Th:BigRigidity} are satisfied, hence the rigidity of the lattice.

  Finally, for all $n \neq 3$ the lattice $\Gamma$ acts transitively on
  $(\Gamma,S)$ thus, by \cref{Prop:LGRSLGR}, it is SLG-rigid.
\end{proof}

%%%%%%%%%%%%%%%%%%%%%%%%%%%%%%%%%%%%%%%%%%%%%%%%%%%%%%%%%%%%%%%%%%%%%%%%%%%%%%%%%%%%%% 
\section{Conclusion and open problems}
% \markboth{Conclusion and open problems}{Conclusion and open problems}
Our main result is proved for graphs quasi-isometric to the Bruhat-Tits building
of $PSL_n(\bK)$ and the key idea of the proof is to use the rigidity of this
building to “transfer it” to the graph quasi-isometric thereto. One can ask
wether we can generalize this idea to other LG-rigid graphs.
\begin{Q}
  Let $\calG$ be quasi-isometric to a LG-rigid graph $\mathcal{H}$, both having
  cocompact isometry group. If the quasi-isometry is
  $\Is(\calG)$-equivariant, is $\calG$ LG-rigid?
\end{Q}
Remark that if $\calH$ and $\calG$ are two Cayley graphs of the same group, we
can chose $\calH$ to be LG-rigid and $\calG$ to be non-rigid (see the discussion
below \cref{CEx:torsion} for more details). In that case the hypothesis of the
preceding question are satisfied without $\calG$ being LG-rigid. Thus, more
restrictive hypothesis will be needed to get the rigidity of $\calG$. 
\bigskip

Our result on lattices is proved for $n\neq 3$; when $n=3$ we don’t
know (yet) the answer. Indeed, our proof is based on 
the rigidity of the Bruhat-Tits building of $PSL_n(\bK)$, a result known to
be true only for $n\neq 3$. In the $n=3$ case, a lot of flexibility seems to be
allowed (see for example \cite{Barre}) obstructing any local recognizability
result. Hence the following question:
\begin{Q}
  Are torsion-free lattices of $SL_3(\bK)$ LG-rigid? 
\end{Q}

Lattices in $p$-adic Lie groups can be viewed as particular cases of
$S$-arithmetic lattices.
\begin{Def} Let $S$ be a set of prime.\\
  We say that $\Gamma$ is an \textbf{$S$-arithmetic} lattice if it’s a lattice in a product of
  the form $\prod_{i} G_i$ where $G_i$ is either a real Lie group or a $p$-adic
  Lie group for $p \in S$.
\end{Def}
Hence, one we can ask what happens in that more general case.
\begin{Q}
  Are torsion-free $S$-arithmetic lattices LG-rigid?
\end{Q}

A result by Bader, Furman and Sauer \cite[Theorem B]{BFS} can be
used to deal with irreducible torsion-free $S$-arithmetic lattices. Indeed, if
the product $\prod_i G_i$ contains at least a non-compact real factor, then the
aforementioned theorem implies that the isometry group of a Cayley graph of 
$\Gamma$ is discrete. Thus, by \cref{Th:DlSTdiscrete} the lattice is LG-rigid.
Now, if the product contains a compact real factor then the isometry group of
the Cayley graph might not be discrete and in that case, the problem is still open.

When the lattice is reducible, we now know that the projection on the $p$-adic factors
gives LG-rigid lattices. Moreover, if we suppose the real factors to be simple
and connected, then a result by de la Salle and Tessera 
\cite{DlST2015} shows that the projection on these factors are also LG-rigid. Hence it remains to
understand how to combine these results on the \emph{factors} in order to get a
result on the \emph{product}.

% \addcontentsline{toc}{section}{Conclusion}
\newpage
\bibliographystyle{alpha}
\bibliography{LGR}
\newpage

\section*{Notations Index}
\markboth{Notations}{Notations}
\addcontentsline{toc}{section}{Notations}
\begin{description}
\item[$\mathfrak{A}$] Atlas of isometries from $\mY$ to $\myg$.
\item[$\Appart$] An apartment in $\bdg$.
\item[$(\Gamma,S)$] Cayley graph of $\Gamma$ with respect to the generating
  part $S$.
\item[$H_0$] The group $\rho^{-1}(PSL_n(\bK))$.
\item[$\Is(\calG)$] Isometry group of $\calG$.
\item[$\iota_y$] Isometry from $\CalY$ to $\bdg$ based at $y$ (see \cpageref{firstmentionofSLGX}).
\item[$\kappa$] Natural injection of $\mY$ in $\CalY$ (see \cref{sec:isommYmyg}).
\item[${[L]}$] Class modulo homothety of the lattice $L$.
\item[$\Print (x)$] The print of the vertex $x$ (see \cref{Def:Print}).
\item[$\PY$] A print in $\mY$ (see \cref{Def:PrintY}).
\item[$\phi_y$] Local isometry from $\CalY$ to $\bdg$ based at $y$ (see \cref{eq:locisomdeCalYversX0}).
\item[$q$] Quasi-isometry between $\myg$ and $\bdg$.
\item[$R$] Radius such that $\mY$ is $R$-locally the same as $\myg$.
\item[$\rho$] Injective homomorphism from $\Is(\myg)$ to $\Is(\bdg)$.
\item[$\ratlas$] See \cref{Prop:ConstructionDeLAtlas}.
\item[$\rprint$] Radius considered to define prints (see \cref{Def:PrintY}).
\item[$\Rbdg$] Radius such that $\CalY$ is $\Rbdg$-locally $\bdg$.
\item[$\rX$] Radius such that $\iota_y$ coincide with $\phi_y$ on
  $B_{\CalY}(y,\rX)$ (see \cpageref{firstmentionofSLGX}).
\item[$\rY$] See \cref{Claim:localeisometrie}. 
\item[$\tau(x)$] The type of the vertex $x$, where $x$ belongs to the
  Bruhat-Tits building of $PSL_n(\bK)$.
\item[$\bdg$] The Bruhat-Tits building of $PSL_n(\bK)$. 
\item[$\CalY$] Hybrid graph built to be locally the same as the building
  (see \cref{Sec:buildingY}).  
\item[$(y_1, \ldots, y_l)$] A path of adjacent vertices $y_1$, $y_2$, …, $y_l$.
\end{description}
\vfill

\noindent\textbf{Amandine Escalier}\\
Université Paris Cité and Sorbonne Université,\\
CNRS, IMJ-PRG,\\
F-75013 Paris\\
France
\end{document}